\documentclass[onecolumn]{article}
\usepackage{algorithm, algcompatible}
\usepackage{graphicx} 
\usepackage{amsmath} 
\usepackage{amssymb}  
\usepackage{amsthm}
\usepackage{amsfonts}
\usepackage{dsfont}
\usepackage{color}
\usepackage{subfig}
\usepackage{pgfplots}
\usepackage{tikz}
\usetikzlibrary{arrows}
\usetikzlibrary{intersections}
\usepackage[english]{babel}
\usepackage[latin1]{inputenc}

\newcommand{\soft}{\mathbb{S}}
\newcommand{\G}{\mathcal{G}}
\newcommand{\V}{\mathcal{V}}
\newcommand{\EE}{\mathcal{E}}
\newcommand{\cardV}{|\mathcal{V}|}
\newcommand{\X}{\mathcal{X}}

\newcommand{\argmin}[1]{\underset{#1}{\mathrm{argmin\,}}}

\newcommand{\xmin}{x^{\star}}
\newcommand{\xmintt}{x_t^{\star T}}
\newcommand{\Xmin}{X^{\star}}

\newcommand{\xtrue}{\widetilde{x}}
\newcommand{\zmin}{z^{\star}}

\newcommand{\ix}{\mathring{x}} 
\newcommand{\iz}{\mathring{z}} 
\newcommand{\iu}{\mathring{u}} 
\newcommand{\reg}{\mathrm{\mathbf{Reg}}^d_T}
\newcommand{\regs}{\mathrm{\mathbf{Reg}}^s_T}
\newcommand{\R}{\mathbb{R}}

\newcommand{\prox}{\mathrm{prox}}
\newcommand{\Ms}{M^{\star}}

\newcommand{\Mq}{M_Q}
\newcommand{\Mp}{M_{\phi}}
\newcommand{\Nv}{\mathcal{N}_v}
\newcommand{\Nw}{\mathcal{N}_w}

\newtheorem{lemma}{Lemma}
\newtheorem{proposition}{Proposition}
\newtheorem{theorem}{Theorem}
\newtheorem{corollary}{Corollary}

\newtheorem{assumption}{Assumption}

\title{\Large \bf  Centralized and distributed online learning \\ for sparse time-varying optimization}

\author{Sophie M. Fosson, \emph{Member, IEEE} }

\begin{document}
\maketitle
\begin{abstract}
The development of online algorithms to track time-varying systems has drawn a lot of attention in the last years, in particular in the framework of online convex optimization. Meanwhile, sparse time-varying optimization has emerged as a powerful tool to deal with widespread applications, ranging from dynamic compressed sensing to parsimonious system identification. In most of the literature on sparse time-varying problems, some prior information on the system's evolution is assumed to be available. In contrast, in this paper, we propose an online learning approach, which does not employ a given model and is suitable for adversarial frameworks. Specifically, we develop centralized and distributed algorithms, and we theoretically analyze them in terms of dynamic regret, in an online learning perspective. Further, we propose numerical experiments that illustrate their practical effectiveness.
\end{abstract}
%
\section{Introduction}
Time-varying optimization has attracted an increasing attention in the last years in machine learning, control, and signal processing, motivated by the observation that usually real-world systems vary with time. Examples are widespread, including big data streams \cite{hal15}, model predictive control \cite{jer14}, resource allocation \cite{che17}, online learning \cite{mok16}, dynamic identification \cite{fox18cdc}, and tracking moving agents \cite{bed18}. Other applications are illustrated in \cite{sim18acc}, and in the recent survey \cite{dall19}.

Formally, by time-varying optimization, we mean a sequence of optimization problems of the kind $\min f_t $, where $t=0,\dots,T$ is the time variable. If the problem can be solved off-line, i.e., after time $T$, then it can be considered as static. Usually, this is not the case: the goal is to track the optimal points as long as the optimization problem varies. This calls for online algorithms, that provide solutions in the system's time-scale, which might be very fast. Moreover, the minimization of $f_t$ might involve the processing of large data; therefore, the development of prompt tracking strategies is challenging.

The literature on online algorithms for time-varying systems is mainly settled in convex optimization, which encompasses a number of applications and is mathematically tractable. To mention some examples, convex functionals are used to model problems of online system identification in \cite{bra16,pra16}, tracking moving targets in \cite{linrib14}, and dynamic magnetic resonance imaging in \cite{vas16}. Most of the theoretical analyses on time-varying convex optimization are oriented towards evaluating the tracking error at time $t$. If $\xmin_t$ is the minimizer of $f_t$ and $x_t$ is the estimate provided by the online algorithm, the tracking error can be defined as a distance between $\xmin_t$ and $x_t$, see, e.g., \cite{linrib14,sim16TSP,rah17,sun17}. As the system is time-varying, the tracking error is not expected to converge to zero in time: an algorithm is considered successful if it guarantees a bounded tracking error, with a sufficiently small bound. 

A different research line addressed, e.g., in \cite{hal15,mok16,che17,sha18,fox18cdc}, proposes the analysis of the \emph{dynamic regret}, a popular performance metric in \emph{online learning}, see \cite{sha12book}. Online learning is the sub-field of machine learning that aims at iteratively learning time-varying models, by assuming a game theoretic or adversarial framework. Specifically, a sequence of rounds is considered: at each iteration $t$, a learner plays an action $x_{t}$; then, an adversary chooses and reveals a loss function $f_t$. Thus, the learner suffers a loss $f_t(x_t)$. In turn, given the knowledge of $f_t$, the learner plays a new action $x_{t+1}$, and suffers a new loss $f_{t+1}(x_{t+1})$, and so on. The dynamic regret is the cumulative difference between the learner's loss $f_t(x_t)$ and $f_t(\xmin_t)$, the last one being the best possible loss in hindsight. If compared to the tracking error metric, the dynamic regret provides more insight on the tracking properties of an online algorithm with respect to the system's evolution.  For example, let us assume to know that the tracking error is bounded by a constant $b$, which is the result of most analyses, as mentioned above. Then, let us assume that the system, after a transient period, tends to converge to a constant value. In this case, a good online algorithm is expected to converge as well. Nevertheless, the knowledge of $b$ does not capture a possible convergent behavior. In contrast, the dynamic regret captures this feature; see, e.g., \cite{mok16} for more details. For this motivation, analyses of the dynamic regret have been gaining an increasing attention in the literature on online learning and online convex optimization, see, e.g, the recent works \cite{zha17,zha18,lu19,dix19}.

The online algorithms proposed in the literature for time-varying convex optimization are usually based on iterative procedures. Both centralized and distributed approaches are developed. Concerning the centralized methods, in \cite{hal15}, a dynamic mirror descent is proposed; in \cite{mok16}, a gradient descent strategy is analyzed; in \cite{che17}, a modified online saddle-point scheme is studied; in \cite{sim17tsp,sim19},  prediction-correction methods are developed for constrained problems. At the same time, several distributed schemes are proposed:  in \cite{linrib14}, decentralization is based on the alternating direction method of multipliers (ADMM); in \cite{sim17}, prediction-correction methods are extended to networked systems; in \cite{sun17}, distributed gradient-based methods are proposed for quadratic problems; in \cite{rah17}, the focus is on continuous-time models; in \cite{lee18}, two decentralized variants of Nesterov primal-dual algorithm are developed; in \cite{sha18}, the mirror descent strategy is decentralized. In \cite{lu19}, a distributed approach based on auxiliary optimization is provided.
%

Within time-varying convex optimization, an important subset is represented by \emph{sparsity} promoting problems, that is, problems whose solution is induced to be a vector with many zero entries. 
Sparsity is nowadays widely studied as it makes it possible to build parsimonious models from large data. In system identification, machine learning, the call for parsimonious models is rapidly increasing to deal with the increasing complexity of systems or with the need of running in small devices, like smartphones.
In signal processing, sparse convex optimization has gained a lot of attention with the advent of compressed sensing (CS, \cite{don06}), which states that  sparse signals can be recovered from few linear measurements. In system identification, the CS paradigm is exploited in the estimation of sparse ARX models from a limited number of observations, see \cite{fox18cdc,tot11,san11}. 

The literature on sparse time-varying optimization (STVO) is quite recent and mainly focused on dynamic CS. Most of the works on the topic assume some prior information on the system's evolution. In \cite{zac12,zin13,mot15,cha16}, the aim is to track time-varying sparse signals which evolve according to Markov models; a Kalman filtering approach is exploited. In \cite{bal15}, a finite bound for the tracking error is assessed, under boundedness assumptions on the  the signal and its derivative. We refer the reader to \cite{vas16} for a complete review.

The goal of this paper is to develop novel strategies and theoretical results for STVO, in terms of dynamic regret, without  prior information on the dynamics. The lack of prior information can be interpreted as an adversarial framework, where the functional is modified arbitrarily. However, it is intuitive that a completely disordered evolution cannot be tracked: the online estimation performance is expected to improve in case of slowly varying systems.

In the literature on time-varying convex optimization, most of the theoretical results exploit the assumption that the $f_t$'s are differentiable, see, e.g., \cite{hal15,mok16,sim17,sim19,sha18,lu19}. This can not be applied to sparse problems, which usually envisage an $\ell_1$ regularizer. Few results do not require differentiability \cite{ber18}, and are limited to prove the boundedness of the tracking error. To the best of our knowledge, no theoretical dynamic regret analysis has been yet performed on non-differentiable functionals. 

In summary, in this paper we propose an approach that differs from most previous literature because, firstly, it provides a dynamic regret analysis that does not require differentiability, and, secondly, it does not assume a specific evolution structure. A first work that combines these two features is \cite{fox18cdc}, where an online algorithm based on iterative soft thresholding (IST) is proposed for STVO from compressed measurements and analyzed in terms of dynamic regret. This paper extends \cite{fox18cdc} in two main directions. First, we propose and analyze a different centralized online algorithm, based on  the Douglas-Rachford splitting \cite{lio79,boy10,gis17}. Second, we develop and analyze an online distributed algorithm, which extends the distributed iterative soft thresholding algorithm (DISTA) proposed in \cite{rfm15} to the time-varying setting. For both  centralized and distributed algorithms, we study the dynamic regret and we present numerical simulations to illustrate their practical effectiveness.

The paper is organized as follows. In Section \ref{sec:ps}, we introduce our specific formulation of STVO. In Section \ref{sec:dr}, we illustrate the online strategy based on splitting, and in Section \ref{sec:dr_dr}, we theoretically analyze it terms of 
dynamic regret. In Section \ref{sec:dista}, we propose the distributed online strategy, which is analyzed in Section \ref{sec:dr_dista}. Section \ref{sec:nr} is devoted to numerical simulations. Finally, we draw some conclusions.
\section{Problem statement}\label{sec:ps}
In this paper, we consider STVO problems that can be modeled as follows:
\begin{equation}\label{newquad}
\begin{split}
\min_{x\in\R^n}&f_t(x),~~~t=0,\dots,T\\
&f_t(x):=h_t(x)+\|x\|_1
\end{split}
\end{equation}
where $h_t$'s are quadratic and strongly convex:
\begin{equation}\label{hquad}
h_t(x):=\frac{1}{2}x^TQ_tx+\phi_t^Tx+c_t
\end{equation}
where $Q_t\in\R^{n,n}$ is symmetric positive definite, $\phi_t\in\R^n$ and $c_t\in\R$. Since $c_t$ is not relevant in the minimization, we neglect it in the rest of the paper. 

Strong convexity is often exploited for time-varying optimization, because it implies contractivity, hence stronger convergence properties, as studied in \cite{mok16,linrib14,hal15,lee18,sim17,sim19,sha18}. In line with these works, we assume strong convexity. Nevertheless, we mention that in \cite{sim17pre}, a possible alternative to strong convexity is illustrated, based on bounded $\alpha$-averaged operators, which might be investigated in future work.

Problem \eqref{newquad}-\eqref{hquad} is the basis for a large class of STVO problems; we illustrate some examples.
\subsubsection{Elastic-net}
Let us consider a time-varying CS problem: given $t=0,\dots, T$, we aim at the online recovery of a sparse signal $\xtrue_t\in\R^n$ (i.e., $\xtrue_t$ has $k_t\ll n$ non-zero components) from compressed, linear measurements. More precisely, at each $t$, we observe
\begin{equation}\label{acquisition}
y_t=A_t \xtrue_t+e_t,~~A_t\in\R^{m, n},~m<n,
 \end{equation}
 where $e_t\in\R^m$ is a possible measurement noise. The goal is to recover $\xtrue_t$ given $y_t$ and $A_t$, knowing that $\xtrue_t$ is sparse. As illustrated in \cite{fox18cdc}, this model can be applied for compressed system identification of linear systems with time-varying parameters. We specify that the exact knowledge of $k_t$ is not required. According to CS theory, an efficient way to tackle this problem is the convex relaxation called  Lasso, which consists in the minimization of $\frac{1}{2}\left\|y_t-A_t x\right\|_2^2+\lambda\left\|x\right\|_1$, where $\lambda>0$; see \cite{tib96,for10} for details. The presence of the $\ell_1$-norm regularizer supports sparsity. If the number of measurements $m$ is sufficiently large, $\xtrue_t$ can be recovered via Lasso, with a bias proportional to the design parameter $\lambda$. We refer the reader to \cite{fou13} for a complete overview on this topic.
 
As $m<n$, the least-squares term in Lasso is not strongly convex. A variation of Lasso, known as Elastic-net \cite{zou05}, enjoys this property by the addition of a Tikhonov $\ell_2$-norm regularizer:
\begin{equation}\label{eq:elasticnet}
	\begin{split}
	&f_t(x):=\frac{1}{2}\left\|y_t-A_t x\right\|_2^2+\frac{\mu}{2}\left\|x\right\|_2^2+\lambda\left\|x\right\|_1\\
	&\lambda>0, \mu>0\\
	\end{split}
\end{equation}
As a difference from Lasso, Elastic-net  promotes a grouping effect of correlated variables instead of selecting just one of them and discarding the others. Moreover, the solution of Lasso necessarily has no more than $m$ non-zero values, see \cite{tib13}, while this limitation is not present in Elastic-net. The effectiveness of Elastic-net is exploited in many applications, ranging from micro-array classification \cite{zou05} to indoor localization \cite{kha17}.
We remark that $\lambda$ and $\mu$ might be time-varying as well. In particular, $\lambda$ is generally designed by using prior knowledge on the sparsity level $k_t$: the higher the sparsity is, the higher  the weight of the $\ell_1$ should be, see, e.g., \cite{fou13}. In this paper, we assume for simplicity that these parameters are constant.
\subsubsection{MPC with sparse control}
Quadratic, strongly convex models are usually exploited in MPC to predict a dynamic system behavior and optimize its control. Recently, the problem of reducing the number of active control inputs in MPC has been gaining an increasing interest in the literature, see, e.g., \cite{gal16,agu17}. This is known as sparse control, and can be tackled by introducing an $\ell_1$ regularizer. The final aim of sparse control  is to reduce consumption and transmission costs.
\subsubsection{Sparse iterative learning control}
The problem of reducing the number of control inputs is investigated in iterative learning control (ILC) as well. The purpose of ILC is the online optimization of repeated systems, with outstanding application in robotics and mechatronics. As in MPC, quadratic cost functionals are widely exploited in ILC \cite{lee00}. In \cite{oom17}, the use of $\ell_1$ regularizers is proven to be effective to obtain a reliable sparse control.
\subsection{Performance metric: dynamic regret}
 Our ultimate goal is to solve Problem \eqref{newquad}-\eqref{hquad}. More precisely, we are interested in computing the minimizer $\xmin_t=\argmin{x\in\R^n} f_t(x)$, which  represents to variable to track. In principle, the problem can be solved at each $t$ through any convex optimization method, as the $f_t$'s are convex. However, we aim to solve the problem online, that is, $\xmin_t$ should be estimated between instant $t$, when $y_t$ and $A_t$ are revealed, and instant $t+1$, when the next data acquisition is performed. Therefore, running a convex optimization algorithm might be not feasible if the time-scale is fast and the dimension $n$ is large. We thus aim at developing fast, suboptimal strategies to track the minima with satisfactory accuracy. The first step to pursue this goal is to choose a suitable performance metric.

In game theory and online learning, a popular performance metric is the dynamic regret, denoted by $\reg$, which is defined as follows (see, e.g., \cite{zin03,mok16}):
\begin{equation*}
\reg (\xmin_1,\dots,\xmin_T):=\sum_{t=1}^{T}\big(f_{t}(x_{t})-f_{t}(\xmin_t)\big)
\end{equation*}
where 
\begin{equation}\label{def:zt}
\xmin_t:=\argmin{x\in\X}f_{t}(x)
\end{equation}
and $x_t$ is the action played by the online algorithm in $[t-1,t)$, thus before that $f_t$ is revealed. $\X\subseteq \R^n$ is the feasible space. Intuitively, an online algorithm is successful if its $\reg$ is sublinear, because this implies that, on average, it performs as well as the clairvoyant opponent that plays the optimal action $\xmin_t$ \cite{haz07,clairvoyant}. The possibility of achieving a sublinear $\reg$ depends on the evolution and regularity of the $f_t$'s. A quantity that well captures the system's evolution is the \emph{path length}, defined as the cumulative distance between reference points \cite{mok16,clairvoyant}; in our setting, we can consider  $\sum_{t=1}^T\left\|\xmin_t-\xmin_{t-1}\right\|_2$ as path length.

As to online algorithms, the following requirements are fundamental: (a) if, for each $t$, the distance $\left\|\xmin_t-\xmin_{t-1}\right\|_2$ between successive minima is bounded, then the estimation error is bounded; (b) if $\left\|\xmin_t-\xmin_{t-1}\right\|_2$ tends to zero, i.e., the system tends to converge, also the estimation error should tend to zero. Requirements (a) and (b) are well captured by the dynamic regret, as illustrated, e.g., in \cite{zin03,mok16}. 
 \subsection{Summary of previous literature on regret analysis}
 Before introducing the algorithms, we briefly overview the previous results on regret analysis in time-varying optimization and online learning. 
 \begin{table}[ht]
\centering
\caption{Main results on static and dynamic regret in the  literature. Each row of the table represents a paper. We distinguish convex (C) and strongly convex (SC) models, and we indicate if algorithms are gradient-based (G) or else. Finally, we specify the main assumptions: $\beta>0$ is a suitable fixed value; $\eta_t\to 0$ denotes the need for vanishing learning parameters, which is undesired for large or infinite time horizons; $\Delta_t:=\left\|x^{\star}_t-x^{\star}_{t-1}\right\|_2$.} 
\resizebox{\columnwidth}{!}{%
\begin{tabular}{|c|c|c|c|c|l|}
 \hline 
  & $f_{t}$ & Alg. &$\regs$ & $\reg$ & Assumptions\\
 \hline  
 \cite{zin03} & C & G &$O(\sqrt{T})$ & $O\left(\sqrt{T}(1+\sum \Delta_t)\right)$ & $\|\nabla f_{t}\|\leq \beta$, $\eta_t\to 0$\\
 \hline  
 \cite{haz07} & SC &Newton &$O(\log T)$ & & $\|\nabla f_{t}\|\leq \beta$, $\eta_t\to 0$\\
 \hline 
\cite{duc10}  & C; SC &COMID& $O(\sqrt{T}); O(\log T)$ & & $\|\nabla h_{t}\|\leq \beta$  \\
 \hline	 
\cite{wan12} & C; SC &ADMM& $O(\sqrt{T}); O(\log T)$ & & $\|\nabla h_{t}\|\leq \beta$\\
 \hline 
 \cite{suz13} & C; SC &ADMM& $O(\sqrt{T}); O(\log T)$ & & $\|\nabla h_{t}\|\leq \beta$, $\eta_t\to 0$\\
 \hline 
 \cite{hos16} & C & ADMM & $O(\sqrt{T})$ & & $\|\nabla h_{t}\|\leq \beta$, $\eta_t\to 0$\\
 \hline 
 \cite{mok16} & SC &G &  & $O(1+\sum \Delta_t)$ & $\|\nabla f_{t}\|\leq \beta$\\
 \hline 
 \cite{fox18cdc} & El.-net (SC) & O-IST  & & $O\left(1+\sum \Delta_t + \sum \Delta_t^2\right)$ &  \\
 \hline
 \cite{akb19} & C & ADMM & $O(\sqrt{T})$ & &  \\
 \hline
 \end{tabular}}\label{table1}
\end{table}

 Minimization problems of the kind $\min_{x\in\R^n}\sum_t f_t(x)$, with $f_t(x)=h_{t}(x)+r(x)$ are widely considered in the literature, where $h_t$ and $r$ are convex, and $r$ is a static regularizer. These problems are intrinsically static: new indirect data are acquired at each time $t$ to estimate a static  optimization variable. This is usually  analyzed in terms of \emph{static} regret, which is defined as  $\regs := \min_x\sum_{t=1}^T\big( f_t(x_t)-f_t(x)\big)$. 

In Table \ref{table1}, we summarize the main results (in chronological order) on regret analysis in the literature. Even though our interest is in the dynamic regret, we also report results on static regret for completeness. 

The main algorithms proposed for the static problem are COMID  \cite{duc10}, based on mirror descent, and different variants of online ADMM in \cite{wan12,suz13,hos16}; in particular, in \cite{hos16} a distributed setting is considered. As shown in Table \ref{table1}, all these methods achieve a static regret of order $O(\sqrt{T})$ for convex cost functionals, and in some cases an improvement to $O(\log{T})$ is obtained in case of strong convexity. As illustrated in the table, most of these methods are driven by decreasing sequences of parameters, which are generally exploited to improve the convergence properties. However, this tool can not be used for tracking problems, where $T$ is possibly infinite.
 
Less work is devoted to the dynamic regret analysis in tracking problems. The main result is provided by \cite{mok16}, which analyzes the dynamic regret of a gradient descent method for strongly convex functionals. In \cite{fox18cdc}, an online iterative soft thresholding (O-IST) method obtains similar performance limited to the Elastic-net model, while not requiring a bounded subgradient.
\subsection{O-IST for quadratic problems}
Before presenting the main algorithms, we retrieve O-IST proposed in \cite{fox18cdc} for Problem \eqref{eq:elasticnet},  and we generalize it to Problem \eqref{newquad}-\eqref{hquad}. This paragraph provides the background to understand the distributed algorithm presented in Section \ref{sec:dista}.

O-IST consists in performing a soft thresholding iteration at each $t$. This is a successful strategy in the sense that $\reg$ (see Table \ref{table1}) is sublinear whenever the path length is sublinear. In particular, this implies that (a) $x_t$ converges to $\xmin_t$ when $\left\|\xmin_t-\xmin_{t-1}\right\|_2$ is null or decreasing as, e.g., $\frac{1}{t^{\beta}}$, $\beta\in(0,1]$, and (b) we have a steady state tracking error for $\left\|\xmin_t-\xmin_{t-1}\right\|_2>c$, for some $c>0$.

In Algorithm \ref{alg:IST_dynamic}, we adapt O-IST \cite{fox18cdc} to Problem \eqref{newquad}-\eqref{hquad}. Differently from \cite{fox18cdc}, we run  $r\geq 1$ soft thresholding iterations at each $t$, where $r$ depends on the time available between $t$ and $t+1$. In Algorithm \ref{alg:IST_dynamic}, the operator $\soft_{\beta}:\R^{n}\to\R^{n}$, $\beta>0$, is the component-wise soft thresholding operator, defined as follows: for $z\in\R$, $\soft_{\beta}[z]=z-\beta$ if $z>\beta$; $\soft_{\beta}[z]=z+\beta$ if $z<-\beta$; $\soft_{\beta}[z]=0$ otherwise; see, e.g., \cite{for10} for details. Moreover, the dynamic regret analysis for O-IST for Problem \eqref{newquad}  can be straightforward derived from the results in \cite{fox18cdc}.

Step 4 of Algorithm \ref{alg:IST_dynamic} is derived as follows. Given a generic problem $\frac{1}{2}x^TQx+\phi^Tx+\lambda\|x\|_1$, a direct minimization over $x$ is not possible, due to the presence of both the $\ell_1$ term and the term $x^TQx$ which couples the variables. In order to decouple the variables,  a surrogate term $\frac{1}{2} (x-b)^T\left[\frac{1}{\tau}I-Q\right](x-b)$ can be added, where $b\in\R^n$ is an auxiliary variable and $\tau$ is designed such that $\frac{1}{\tau}I-Q$ is positive definite. In this way, the surrogate term is always non negative, and the global minimum is the same by adding it. An alternated minimization is then performed: with respect to $b$, the minimum is obtained for $b=x$; with respect to $x$, the problem is decoupled and can be solved by soft thresholding \cite{for10,fox18cdc}.
\begin{algorithm}
	\caption{O-IST for Problem \eqref{newquad}-\eqref{hquad}}
	\label{alg:IST_dynamic}
	\begin{algorithmic}[1] 
		\STATEx {\bf{input}}: $\lambda>0$, $\tau>0$, $x_0=0$; at time $t=0,\dots,T$, $Q_t$ and $\phi_t$;
		\STATEx {\bf{output}}: in $[t,t+1)$, an estimate $x_{t+1}$ of $\xmin_t$;
		\FOR{$t=0,\dots,T$}
		\STATE $\ix_0=x_{t}$; $Q_t$ and $\phi_t$ are revealed;
		\FOR{$h=1,\dots,r$}
		%
		\STATE $\ix_{h}=\soft_{\lambda}\left[\ix_{h-1}-\tau Q_t\ix_{h-1}-\tau\phi_t\right]$
		\ENDFOR
		\STATE $x_{t+1}=\ix_{r}$
		\ENDFOR
	\end{algorithmic}
\end{algorithm}
O-IST is successful in terms of dynamic regret; however, in practice IST methods are observed to be not very fast, which, in the online version, reduce the promptness to sudden changes  \cite{fox18cdc}. For this motivation, in this paper we develop a faster online strategy based on Douglas-Rachford splitting, whose convergence properties in the static framework can be leveraged to obtain good tracking properties in the dynamic framework. 

We remark that accelerated versions of IST might be investigated as well to speed up O-IST, based, e.g., on Nesterov accelerations \cite{sop16} or FISTA \cite{bec09}. However, these methods are driven by time-varying, convergent parameters, which makes their application more difficult in a tracking context,  where the time horizon is possibly infinite. For this motivation, we focus on splitting methods.

Concerning the distributed setting, in static sparse recovery, a decentralization of IST is proposed  in \cite{rfm15}. By leveraging  \cite{rfm15}, in the second part of this work, we develop a distributed online version of IST and we analyze its dynamic regret.

We specify that an online distributed splitting methods could be conceived as well, as distributed/parallel splitting algorithms are widely applied in sparse optimization, see, e.g., \cite{mata15,fia18}. However, a rigorous dynamic regret analysis of an online distributed splitting is rather technical, thus left for future work.
\section{O-DR: Online Douglas-Rachford splitting}\label{sec:dr}
In this section, we present an online splitting algorithm to tackle Problem \eqref{newquad}, based on the Douglas-Rachford (DR) method \cite{lio79,gis17}. First, we briefly review the classical batch DR  algorithm in a static framework. 

DR is an iterative algorithm that tackles the minimization of cost functionals of the kind $f(x)=h(x)+g(x)$, $x\in\R^n$. 
The procedure can be formulated as follows. Given the proximal operator, defined by
\begin{equation}\label{proximal}
 \prox_{\gamma h }(z):=\argmin{x\in\R^n}\left[\gamma h(x)+\frac{1}{2}\|x-z\|_2^2\right],
\end{equation}
where $\gamma>0$, for each $t=0,\dots,T_{stop}$,
\begin{equation}\label{dr_procedure}
 \begin{split}
u_{t} & =\prox_{\gamma g} (2x_{t}-z_{t})\\
z_{t+1} & =z_{t}+2\alpha(u_{t}-x_{t})\\
x_{t+1} & =\prox_{\gamma h}(z_{t+1})
\end{split}
\end{equation}
where $\alpha>0$, and $T_{stop}$ is the instant where some stop criterion is met. 
The procedure can be equivalently written as
\begin{equation}\label{dr_procedure_short}
 \begin{split}
z_{t+1} & =\mathcal{R}(z_{t})\\
x_{t+1} & =\prox_{\gamma h}(z_{t+1})
\end{split}
\end{equation}
where $\mathcal{R}:=(1-\alpha)I+\alpha(2\prox_{\gamma h}-I)(2\prox_{\gamma g}-I)$, $I$ being the identity operator.
The sequence $x_t$ is proven to converge to $\xmin=\argmin h(x)+g(x)$ (while $z_t$ converges to a fixed point of $\mathcal{R}$: $\zmin=\mathcal{R}\zmin$) when $h$ and $g$ are proper closed and convex, and $\alpha\in(0,1)$, see \cite{gis17} and references therein for more details. The case $\alpha=1$, also known as Peaceman-Rachford splitting method, converges faster than the case $\alpha\in(0,1)$, under strong convexity assumptions on $h$, see \cite[Section III]{gis17}. 

We remark that DR is equivalent to ADMM \cite{boy10} for convex problems $h(x)+g(x)$. More precisely, ADMM tackles more general problems of kind $h(x)+g(z)$ subject to $Ax+Bz=c$, and actually is the dual version of DR. The equivalence between DR and ADMM is widely studied in the literature, see, e.g., \cite{mou19} and references therein. In this paper, we leverage the DR formulation; however, the ADMM formulation is possible as well, and yields the same theoretical and numerical results.

To apply DR to Problem \eqref{newquad} (for the moment, in the static case $Q_t=Q$, $\phi_t=\phi$), we set $h(x)=x^T Q x+\phi^T x$, and $g(x)=\lambda\|x\|_1$. The steps in \eqref{dr_procedure} are made explicit in Algorithm \ref{alg:batch_DR}. To unburden the notation, we set $\gamma=1$, and we consider the Peaceman-Rachford version $\alpha=1$ \cite{gis17}. These values are observed to be suitable for the proposed setting; an optimal tuning is beyond the scope of the paper and left for future analysis. 
\begin{algorithm}
\caption{Batch DR for Problem \eqref{newquad}-\eqref{hquad}}
\label{alg:batch_DR}

\begin{algorithmic}[1] 

\STATEx {\bf{input}}: $\lambda>0$, $\mu>0$, $Q\in\R^{n,n}$, $\phi\in\R^n$, $z_0=0$, $x_0=0$

\STATEx {\bf{output}}: at time $T_{stop}$, an estimate $x_{T_{stop}}$ of $\xmin$ 

\FOR{$t=0,\dots,T_{stop}$}


\STATE $u_{t}=\soft_{\lambda}\left[2x_{t}-z_{t}\right]$

\STATE $z_{t+1}=z_{t}+2(u_{t}-x_{t})$

\STATE $x_{t+1}=[Q+I]^{-1}\left[z_{t+1}-\phi\right]$
\ENDFOR

\end{algorithmic}
\end{algorithm}

Afterwards, following the rationale of O-IST, we propose O-DR, that performs $r$ DR steps at each $t$. This is summarized in Algorithm \ref{alg:online_DR}.

\begin{algorithm}
\caption{O-DR for Problem \eqref{newquad}-\eqref{hquad}}
\label{alg:online_DR}

\begin{algorithmic}[1] 
\STATEx {\bf{input}}: $\lambda>0$, $\mu>0$, $z_0=0$; $x_0=0$ at $t=0,\dots,T$: $Q_t$, $\phi_t$
\STATEx {\bf{output}}: in $[t,t+1)$, an estimate $x_{t+1}$ of $\xmin_t$
\FOR{$t=0,\dots,T$}
\STATE $\iz_0=z_{t}$, $\ix_0=x_{t}$, $Q_t$ and $\phi_t$ are revealed
\FOR{$h=1,\dots,r$}
%
%
\STATE $\iu_{h}=\soft_{\lambda}\left[2\ix_{h}-\iz_{h}\right]$
\STATE $\iz_{h+1}=\iz_{h}+2(\iu_{h}-\ix_{h})$
\STATE $\ix_{h+1}=[Q_t+I]^{-1}\left[\iz_{h+1} -\phi_t\right]$
\ENDFOR
%
%
%
\STATE $z_{t+1}=\iz_{r+1}$
\STATE $x_{t+1}=\ix_{r+1}$
\ENDFOR
\end{algorithmic}
\end{algorithm}

\subsection{Related literature}
The research on online splitting methods is very active in these years. In particular, the idea of performing one ADMM/DR iteration at each time step to tackle dynamic problems is known in the literature. However, previous work is mainly focused on the online estimation of a static quantity, which yields to online ADMM/DR procedures different from  O-DR (Algorithm \ref{alg:online_DR}).
Specifically, in \cite{wan12}, an online ADMM, called OADM, is proposed to tackle static problems of kind $\min_x \sum_{t=1}^Tf_t(x)+g(x)$. The idea is to update the estimation of the global minimum at each $t$, when new measurements are acquired, i.e., a new $f_t$ is revealed.  Then, OADM is conceived to tackle the online estimation of a static quantity, which is intrinsically different from the tracking problem proposed in this paper. Similarly to O-DR, OADM perform one ADMM iteration at each time step. Differently from O-DR, in OADM a Bregman divergence term is added to obtain good static regret properties. The analysis in \cite{wan12} is specific for static regret; in particular, some step size parameters of the algorithms are required to increase or decrease as $T$, which can not be applied in a tracking context, where $T$ may be infinite. 
In \cite{suz13}, the same static problem is tackled with slightly different online ADMM procedures, based the addition of proximal operators or $\ell_2$ regularization terms. Similarly to \cite{wan12}, the so-obtained online ADMM procedures work for decreasing/increasing time step parameters, which requires a finite $T$, and prevents their application to tracking problems.
The tracking capabilities of ADMM/DR have been investigated more recently, with particular attention to specific practical problems. In \cite{dall17}, ADMM is used for real-time optimization of power systems; in \cite{mar18}, the tracking capabilities of ADMM are tested  in a dynamic beam-forming problem. The algorithms proposed in these works are based on the idea of performing one ADMM iteration at each time step, and bounds for their limit errors are studied. 
In \cite{mar19}, ADMM is used to track the solution of a stochastic sequence of  problems, parametrized by a discrete time Markov process. Finally, in \cite{cao19}, a dynamic ADMM procedure that performs one ADMM iteration at each time step is analyzed for the dynamic sharing problem: under technical assumptions (in particular,  the time-varying cost functional is sum of strongly convex functions), the convergence to a neighborhood of the optimal time-varying point is proven. Moreover, a numerical experiment on dynamic Lasso is illustrated, even though Lasso does not enjoy the above mentioned technical conditions. 
\section{Dynamic regret analysis for O-DR}\label{sec:dr_dr}
In this section, we show how $\reg$ for O-DR depends on the system's evolution, in terms of path length. This is achieved through some intermediate results.

In this section, we assume that the evolution of the system is bounded.
\begin{assumption}\label{ass:limited}
There exist $\Mq>0$ and $\Mp>0$ such that, for all $t=0,\dots,T$, $\|Q_t\|_2\leq \Mq $ and $\|\phi_t\|_2\leq \Mp$.
\end{assumption}
As a consequence, the minimizer $\xmin_t$ of $f_t$ in Problem \eqref{newquad}-\eqref{hquad} is bounded as well, i.e., there exists $\Ms>0$ such that, for all $t=0,\dots,T$, $\|\xmin_t\|_2\leq \Ms$. 
We remark that this is an assumption on the system's evolution, while we do not force any boundedness on the algorithm's evolution $x_t$ or $f_t(x_t)$. This is an improvement with respect to \cite{mok16}, where the boundedness of $\nabla f_t$ is required \cite[Assumption 3]{mok16}, which excludes, for example, quadratic cost functionals over non-compact state spaces.

In \cite{gis17}, a novel convergence rate analysis is proposed for DR when the cost functional $h(x)+g(x)$ enjoys the additional properties that $h(x)$ is $\sigma$-strongly convex, i.e., $h-\frac{\sigma}{2}\|\cdot\|_2^2$ is convex, and $\beta$-smooth, i.e., $\frac{\beta}{2}\|\cdot\|_2^2-h$ is convex, see \cite[Section II]{gis17} for details. Based on this result, we can prove the following proposition.
 \begin{proposition}\label{prop:contr} Let us consider the static Problem \eqref{newquad}-\eqref{hquad} (that is, $f_t=f$ for each $t$). Let $(x_t,z_t)_{t=0,\dots T_{stop}}$ be the sequence generated by batch DR (Algorithm \ref{alg:batch_DR}) and let $(\xmin,\zmin)$ be its limit point.
Then, $z_t$ converges Q-linearly to $\zmin$; more precisely, for each $t=0,\dots,T_{stop}$
\begin{equation}\label{qlinear}
\left\|z_{t+1}-\zmin\right\|_2\leq \delta \left\|z_{t}-\zmin\right\|_2
\end{equation}
where 
\begin{equation}\label{delta}
\begin{split}
 \delta&=\max \left ( \frac{1-\sigma}{1+\sigma},\frac{\beta-1}{\beta+1}\right)<1 
\end{split}
\end{equation}
where $\sigma$ and $\beta$ respectively are the minimum non-null and the maximum eigenvalues of $Q$.
Moreover,
 \begin{equation}\label{rlinear}
 \left\|x_{t+1}-\xmin\right\|_2\leq q \left\|z_t-\zmin\right\|_2
 \end{equation}
 where $q=\frac{\delta}{1+\sigma}.$
 \end{proposition}
\begin{proof}
The proof is based on \cite[Theorem 2, Corollary 1]{gis17}. First, we notice that \cite[Assumption 2]{gis17} is satisfied by  Problem \eqref{newquad}-\eqref{hquad}: $h(x)$ is $\beta$-smooth and $\sigma$-strongly convex with parameters as defined after \eqref{delta}.
Therefore, \cite[Theorem 2]{gis17} holds, which states that $z_t$ in Algorithm \ref{alg:batch_DR} converges Q-linearly to $\zmin$, with contraction parameter $\delta$ as defined in \eqref{delta}.

Since $\prox_h(z)=(Q+I)^{-1}(z-\phi)$, then $\prox_h$ is $\frac{1}{1+\sigma}$-Lispchitz continuous. Hence, we prove \eqref{rlinear} as follows:
\begin{equation*}
\begin{split}
\left\|x_{t+1}-\xmin\right\|_2&= \left\|\prox_h(z_{t+1})-\prox_h(\zmin)\right\|_2\\
&\leq \frac{1}{1+\sigma}\left\|z_{t+1}-\zmin)\right\|_2\leq  \frac{\delta}{1+\sigma}\left\|z_{t}-\zmin\right\|_2.
\end{split} 
\end{equation*}
\end{proof}
If $r$ iterations are played at time $t$, we easily derive the following corollary from Proposition  \ref{prop:contr}.
\begin{corollary}\label{cor:contr}
For O-DR (Algorithm \ref{alg:online_DR}), the following properties hold:
 \begin{equation}\label{qlineart}
\left\|z_{t+1}-\zmin_t\right\|_2\leq \delta^{r} \left\|z_{t}-\zmin_t\right\|_2
\end{equation}
and
\begin{equation}\label{rlineart}
\left\|x_{t+1}-\xmin_t\right\|_2\leq q^{r} \left\|z_t-\zmin_t\right\|_2. 
\end{equation}
\end{corollary}
By using Corollary \ref{cor:contr}, we  prove that, at  $t$, the distance between the played action and the current minimum is controlled by the distance between successive minima. In the following, we name:
\begin{equation}
\begin{split}
&\Delta_{z_t}:=\|z_t-z^{\star}_{t}\|_2,\\
&\Delta^{\star}_{x_t}:=\|x^{\star}_t-x^{\star}_{t-1}\|_2,~~~\Delta^{\star}_{z_t}:=\|z^{\star}_t-z^{\star}_{t-1}\|_2.
\end{split}
\end{equation}
\begin{lemma}\label{lem:somme}
For O-DR,
	\begin{equation*}
	\begin{split}
	(a)~~&\sum_{t=1}^T \Delta_{z_t}\leq c_1 + c_2 \sum_{t=1}^T \Delta^{\star}_{z_t}\\
	(b)~~&\sum_{t=1}^T \Delta_{z_t}^2\leq c_3 + c_4 \sum_{t=1}^T \Delta^{\star}_{z_t}+c_5 \sum_{t=1}^T \Delta^{\star^2}_{z_t}
	\end{split}
	\end{equation*}
	where 
	\begin{equation}\label{cconstant}
\begin{split}
c_1&= \frac{\delta^r}{1-\delta^r}\left(\Delta_{z_0} - \Delta_{z_T}\right),~~~c_2=\frac{1}{1-\delta^r},\\
c_3&= \frac{\delta^{2r}\left(\Delta_{z_0}^2 - \Delta_{z_T}^2\right)+4 \Ms \delta^{2r}\left(\Delta_{z_0}-\Delta_{z_T}\right)+4\Ms \delta^r c_1}{1-\delta^{2r}},\\
c_4&=\frac{4\Ms \delta^r c_2}{1-\delta^{2r}},~~~c_5=\frac{1}{1-\delta^{2r}}.
\end{split}
\end{equation}
\end{lemma}
\begin{proof}
By the triangle inequality and Corollary \ref{cor:contr}, for each $t=1,\dots,T$,
\begin{equation}\label{ondelta}
\begin{split}
\Delta_{z_t}&=\|z_t-z^{\star}_{t}\pm z^{\star}_{t-1}\|_2\leq \|z_t-z^{\star}_{t-1}\|_2 +\Delta^{\star}_{z_t}\\&\leq \delta^r\Delta_{z_{t-1}} +\Delta^{\star}_{z_t}.
\end{split}
\end{equation}
By summing  over $t=1,\dots,T$, we prove $(a)$:
\begin{equation*}
\begin{split}
&(1-\delta^r)\sum_{t=1}^T \Delta_{z_t}\leq  \delta^r\left(\Delta_{z_0}  -\Delta_{z_T}\right)+ \sum_{t=1}^T \Delta^{\star}_{z_t}.
\end{split}
\end{equation*}
To prove $(b)$, first we use \eqref{ondelta} and the fact that  $\Delta^{\star}_{z_t}\leq \|z^{\star}_t\|_2+\|z^{\star}_{t-1}\|_2\leq 2 \Ms$: 
\begin{equation*}
\begin{split}
\Delta_{z_t}^2& \leq \delta^{2r}\Delta_{z_{t-1}}^2 +\Delta^{\star 2}_{z_t}+2\delta^r\Delta_{z_{t-1}}\Delta^{\star}_{z_t}\\
& \leq \delta^{2r}\Delta_{z_{t-1}}^2 +\Delta^{\star 2}_{z_t}+4\Ms \delta^r\Delta_{z_{t-1}}.
\end{split}
\end{equation*}
Then, we sum over over $t=1,\dots,T$:
\begin{equation}\label{summing...}
\begin{split}
&\sum_{t=1}^T\Delta_{z_t}^2 \leq \hspace{-0.5mm} \frac{\delta^{2r}\left(\Delta_{z_{0}}^2-\Delta_{z_{T}}^2\right)+\sum_{t=1}^T\left(\Delta^{\star 2}_{z_t}+4\Ms \delta^r \Delta_{z_{t-1}}\right)}{1-\delta^{2r}}.
\end{split}
\end{equation}
Since $\sum_{t=1}^T\Delta_{z_{t-1}}=\Delta_{z_{0}}-\Delta_{z_{T}}+\sum_{t=1}^T\Delta_{z_{t}}$, by applying $(a)$, we have
\begin{equation*}
 \sum_{t=1}^T\Delta_{z_{t-1}}\leq \Delta_{z_{0}}-\Delta_{z_{T}}+c_1+c_2 \sum_{t=1}^T \Delta^{\star}_{z_t}.
\end{equation*}
By substituting this bound in \eqref{summing...}, the thesis is obtained, with constants as in \eqref{cconstant}.
\end{proof}
The following two lemmas highlight properties of the quadratic functional in \eqref{newquad}.
\begin{lemma}\label{lem:t_t_1}
For each $t$, and for any $x\in\R^n$, 
\begin{equation}\label{unifcont1}
f_{t}(x)-f_{t}(\xmin_{t})\leq \alpha_1\left\|x-\xmin_{t}\right\|_2+\alpha_2\left\|x-\xmin_{t}\right\|_2^2\end{equation}
where $\alpha_1=\Mq\Ms+\Mp+\lambda\sqrt{n}$ and $\alpha_2=\frac{\Mq}{2}$.
\end{lemma}
\begin{proof} 
We recall that, for any $x,z\in\R^n$ and symmetric $B\in\R^{n,n}$,
\begin{equation}\label{algebra}
x^T B x-z^T B z=(x-z)^T B (x+z).
\end{equation}
By using \eqref{algebra}, we compute the following bound:
\begin{equation}\label{calcolo}
 \begin{split}
&f_{t}(x)-f_{t}(\xmin_{t})=\\
&=\frac{x^T Q_tx-\xmintt Q_t\xmin_{t}}{2}+(x-\xmin_{t})^T\phi_t+\lambda\|x\|_1-\lambda\|\xmin_{t}\|_1\\
&\leq \frac{(x-\xmin_{t})^T Q_t(x+\xmin_{t})}{2}+(x-\xmin_{t})^T\phi_t+\lambda\|x-\xmin_{t}\|_1\\
&\leq \frac{1}{2}\left\|x-\xmin_{t}\right\|_2 \left\|Q_t\right\|_2 \left\|x_t+\xmin_{t}\right\|_2+\left\|x-\xmin_{t}\right\|_2\left\|\phi_t\right\|_2\\
&~~~+\lambda\sqrt{n}\|x-\xmin_{t}\|_2.
\end{split}
\end{equation}
From standard properties of norms, we have
\begin{equation}\label{itsavesme}
\begin{split}
&\|x+\xmin_{t}\|_2\leq \|x-\xmin_{t}\|_2+2\|\xmin_{t}\|_2. 
\end{split}
\end{equation}
Moreover, from Assumption \ref{ass:limited}, we have $\|\xmin_{t}\|_2\leq \Ms$,$\|Q_t\|_2\leq \Mq$, $\|\phi_t\|_2\leq \Mp$.
Therefore, by plugging \eqref{itsavesme} in \eqref{calcolo} and by exploiting Assumption \ref{ass:limited}, we obtain
\begin{equation*}
 \begin{split}
f_{t}(x)-f_{t}(\xmin_{t})\leq& \left\|x-\xmin_{t}\right\|_2 \Mq \left(\frac{\left\|x-\xmin_{t}\right\|_2}{2}+\Ms\right)+\\
&+\left\|x-\xmin_{t}\right\|_2\Mp+\lambda\sqrt{n}\|x-\xmin_{t}\|_2.
\end{split}
\end{equation*}
Then, the thesis is  proven with $\alpha_1=\Mq\Ms+\Mp+\lambda\sqrt{n}$ and $\alpha_2=\frac{\Mq}{2}$.
\end{proof}
Based on these results, we prove the main result.
\begin{theorem}\label{theo:regret} O-DR for Problem \eqref{newquad}-\eqref{hquad} (Algorithm \ref{alg:online_DR}) has the following dynamic regret bound:
\begin{equation*}
\begin{split}
 \reg \leq &\eta_0+\sum_{t=1}^T \left(\eta_1 \Delta^{\star}_{z_t} + \eta_2  \Delta^{\star^2}_{z_t}+\eta_3 \Delta^{\star}_{x_t} + \eta_4 \Delta^{\star^2}_{x_t}\right)
\end{split}
\end{equation*}
where $\eta_i>0$, $i=0,\dots,4$ are assessed in the proof.
\end{theorem}
In particular, this theorem implies that if the path lengths $\sum_{t=1}^T \Delta^{\star}_{z_t}$ and $\sum_{t=1}^T \Delta^{\star}_{x_t}$ are sublinear, then also $\reg$ is sublinear, i.e., the algorithm is successful, in line with previous results shown in Table \ref{table1}.
\begin{proof}
From Lemma \ref{lem:t_t_1}, by considering $x=x_t$,
\begin{equation}\label{nn3}
\begin{split}
f_t(x_t)-&f_t(\xmin_t)\leq \alpha_1\left\|x_t-\xmin_{t}\right\|_2+\alpha_2\left\|x_t-\xmin_{t}\right\|_2^2.
\end{split}
\end{equation}
From Corollary \ref{cor:contr}, we know that $\left\|x_t-\xmin_{t-1}\right\|_2\leq q^r \|z_{t-1}-z^{\star}_{t-1}\|_2=q^r \Delta_{z_{t-1}}$. Then, by applying the triangle inequality, we have
\begin{equation}\label{nn1}
\left\|x_t-\xmin_{t}\right\|_2\leq \left\|x_t-\xmin_{t-1}\right\|_2+\Delta^{\star}_{x_t}\leq q^r \Delta_{z_{t-1}}+\Delta^{\star}_{x_t}.
\end{equation}
Since $(a-b)^2\leq 2a^2+2b^2$ for any $a,b\in\R$, we have
\begin{equation}\label{nn2}
\left\|x_t-\xmin_{t}\right\|_2^2\leq 2q^{2r} \Delta_{z_{t-1}}^2+2\Delta^{\star^2}_{x_t} 
\end{equation}
By substituting \eqref{nn1} and \eqref{nn2} in \eqref{nn3}, we conclude:
\begin{equation*}
\begin{split}
&f_t(x_t)-f_t(\xmin_t)\leq \alpha_1\left( q^r \Delta_{z_{t-1}}+\Delta^{\star}_{x_t}\right) +\\
&~~+ 2\alpha_2 \left(q^{2r} \Delta_{z_{t-1}}^2+\Delta^{\star^2}_{x_t}\right)\\
&=\zeta_1 \Delta_{z_{t-1}}+\zeta_2 \Delta_{z_{t-1}}^2+\zeta_3\Delta^{\star}_{x_t}+\zeta_4\Delta^{\star^2}_{x_t}
\end{split}
\end{equation*}
where 
$\zeta_1= \alpha_1 q^r$, $\zeta_2=2\alpha_2 q^{2r}$, $\zeta_3 = \alpha_1$, $\zeta_6=2\alpha_2.$
%
Now, let us sum over $t=1,\dots,T$:
\begin{equation*}
\begin{split}
\reg&=\sum_{t=1}^T\left(f_t(x_t)-f_t(\xmin_t\right)\\&\leq \kappa + 
\sum_{t=1}^T\left(\zeta_1 \Delta_{z_{t}}+\zeta_2 \Delta_{z_{t}}^2+\zeta_3\Delta^{\star}_{x_t}+\zeta_4\Delta^{\star^2}_{x_t} \right)
\end{split}
\end{equation*}
where $\kappa=\zeta_1(\Delta_{z_0}-\Delta_{z_T})+\zeta_2(\Delta^2_{z_0}-\Delta^2_{z_T})$.
Then, we apply Lemma \ref{lem:somme}:
\begin{equation*}
\begin{split}
&\reg\leq \kappa+ \zeta_1\left(c_1+c_2\sum_{t=1}^T\Delta^{\star}_{z_t}\right)+\zeta_2c_3+\\&+ \zeta_2\left(c_4\sum_{t=1}^T\Delta^{\star}_{z_t}+c_5\sum_{t=1}^T\Delta^{\star^2}_{z_t}\right) +\zeta_3\sum_{t=1}^T\Delta^{\star}_{x_t}+\zeta_4\sum_{t=1}^T\Delta^{\star^2}_{x_t}.
\end{split}
\end{equation*}
Then, the thesis is obtained with $\eta_0=\zeta_1 c_1+ \zeta_2 c_3+\kappa$, $\eta_1=\zeta_1 c_2+\zeta_2 c_4$, $\eta_2= \zeta_2 c_5$, $\eta_3=\zeta_3$, $\eta_4=\zeta_4$.
\end{proof}

%
\section{O-DISTA: Distributed online IST}\label{sec:dista}
As mentioned in the introduction, several works in the literature are concerned with distributed methods to deal with online convex optimization, see, e.g., \cite{linrib14,rah17,sun17,sim17,lee18}. For this motivation, we propose a distributed algorithm for STVO. Specifically, this is a decentralization of the IST algorithm proposed in \cite{fox18cdc}, based on the DISTA algorithm \cite{rfm15}. The key idea is that DISTA can be reformulated for Problem \eqref{newquad}-\eqref{hquad}, and contraction properties can be proven. By starting from this observation, a dynamic regret analysis is performed in Section \ref{sec:dr_dista}.

As in \cite{rfm15}, let us consider an undirected graph $\mathcal{G}=(\mathcal{V, E})$ where $\mathcal{V}$ is the set of nodes, whose cardinality is denoted by $\cardV$, and $\mathcal E\subseteq \mathcal{V\times V}$ is the set of edges. $\mathcal E$ enjoys the property: $(i,j)\in \mathcal E$ implies $(j,i)\in \mathcal E$. $(i,i)\in \mathcal E$ for all $i\in\mathcal V$. A graph is said to be $d$-regular if each node is connected to $d-1$ nodes different from itself. $\mathcal{N}_v$ denotes the neighborhood of $v$: $w\in\mathcal{N}_v$ when $(w,v)\in \EE$.
To lighten the theoretical analysis, we set the following topology conditions.
\begin{assumption}\label{ass:dreg}
 $\G=(\mathcal{V}, \mathcal{E})$ is connected and $d$-regular. Its topology is time-invariant.
\end{assumption}
The connectivity is a natural assumption for collaborative  networked systems. On the other hand,  time-invariance is typical of ideal networks, where no communications interruptions or links failures occur. The $d$-regularity is a technical condition that simplifies the convergence proof, as illustrated in \cite{rfm15,rav15}. In real applications, exact $d$-regularity is uncommon, while an approximated regularity is recommended in many cases, i.e., the presence of clusters or isolated nodes calls for different collaboration protocols. In Section \ref{sub:indoor}, we show an example with non-regular topology. More details on the robustness of thresholding-consensus algorithms to non-regularity and communications losses can be found in \cite{rav15,del16}.


Let $X:=(x_{1},\ldots,x_{|\mathcal{V}|})\in\R^{n,\cardV}$. At each $t$, in the philosophy of \cite[Equation (9)]{rfm15},  we formulate the problem the following problem:
\begin{equation}\label{distributed_t}
\begin{split}
&\min_{X\in\R^{n,\cardV}}F_t(X)\\
&F_t(X):=\sum_{v\in\mathcal{V}}\left[\frac{1}{2}x_v^TQ_{v,t}x_v+\phi_{v,t}^Tx_v+\lambda\|x_v\|_1+\right.\\
&\left.+\frac{1}{2d\tau}\sum_{w\in\mathcal{N}_v}\|\overline{x}_w-{x}_v\|_2^2\right]
\end{split}
\end{equation}
where $\overline{x}_w$ denotes the local mean of $x_w$: $\overline{x}_w=\frac{1}{d}\sum_{w'\in\mathcal{N}_w}x_{w'}$, and $\tau>0$ is a weight that will be assessed later. The  term $\sum_{w\in\mathcal{N}_v}\|\overline{x}_w-{x}_v\|_2^2$ induces a consensus among the local estimates of each node. The motivation to use $\overline{x}_w$ instead of $x_w$ (which would equivalently support the consensus) is rather technical; in a nutshell, it makes easier to split each iteration in two steps: one of local communication and one of individual descent, as illustrated in algorithms \ref{alg:batch_DISTA} and \ref{alg:online_DISTA}.

At time $t$,  each $v\in\V$ is assumed to know local data $Q_{v,t}\in\R^{n,n}$ and $\phi_{v,t}\in\R^{n}$. Nodes aim to track the  minimizer of \eqref{distributed_t}, denoted as $\xmin_t\in\R^n$, by leveraging local information and local communication. 

As for O-IST, the minimization of each $F_t(X)$ can be tackled by considering a surrogate functional, that has the same global minimum of $F_t(X)$, and that can be tackled via alternated minimization, in the presence of local communication. The surrogate functional is obtained from $F_t(x)$ by adding a  quadratic term to deal with the $\ell_1$ term, and by substituting the local mean $\overline{x}_w$ with a local auxiliary variable. More precisely, we define the auxiliary variables $C=(c_1,\dots,c_{\cardV})\in\R^{n,\cardV}$, $B=(b_1,\dots,b_{\cardV})\in\R^{n,\cardV}$, and we define the surrogate functional as follows:
\begin{equation}\label{distributed_t_surro}
\begin{split}
F_t(X,&C,B):=\sum_{v\in\mathcal{V}}\left[\frac{1}{2}x_v^TQ_{v,t}x_v+\phi_{v,t}^Tx_v+\lambda\|x_v\|_1+\right.\\
&\left.+\frac{1}{2d\tau}\sum_{w\in\Nv}\|c_w-{x}_v\|_2^2+\right.\\
&\left.+\frac{1}{2}(x_v-b_v)^T\left(\frac{1}{\tau}I-Q_{v,t}\right)(x_v-b_v)\right].
\end{split}
\end{equation}
As discussed for O-IST, $\tau$ must be designed such that $\frac{1}{\tau}I-Q_{v,t}$ is positive definite. We specify that different $\tau$'s might be considered at each node; here we consider a unique value to simplify the notation.

Now, we can minimize $F_t(X,C,B)$ by alternating minimization over $X\in\R^{n,\cardV}$, $C\in\R^{n,\cardV}$, $B\in\R^{n,\cardV}$. In particular, the minimization with respect to $x_v$ is done by recalling that, given $a>0$, $\argmin{x\in\R}\frac{1}{2}ax^2+bx+\lambda|x|=\soft_{\lambda/a}\left[-b/a\right]$. 

The so-obtained algorithm, which generalizes the method of \cite{rfm15}, is reported in Algorithm \ref{alg:batch_DISTA}. As in the centralized case, we first illustrate the static case $F_t=F$ for each $t$.
\begin{algorithm}  
\caption{Batch DISTA for Problem \eqref{distributed_t}} 
\label{alg:batch_DISTA}
\begin{algorithmic}[1]
\STATEx {\bf{input}}: $\lambda>0$, $\tau>0$; for each $v\in\V$,  $x_{v,0}=0$, $Q_v$;  $\phi_v$
\STATEx {\bf{output}}: at time $T_{stop}$, for each $v\in\V$, an estimate $x_{v,T_{stop}}$ of $\xmin$ 

\FOR{$t=0,\dots,T_{stop}$}
\STATE If $t$ is even, for any $v\in\V$,
     \begin{align*}
 c_{v,t+1}&=\overline{x}_{v,t} \\
x_{v,t+1}&=x_{v,t}
\end{align*}
\STATE If $t$ is odd, for any $v\in\V$,
\begin{align*}
c_{v,t+1}&=c_{v,t}\\
 {x}_{v,t+1}&=\soft_{\lambda\tau/2}\left[\frac{{x}_{v,t}+\overline{c}_{v,t} -\tau Q_{v} {x}_{v,t} -\tau\phi_v}{2} \right]\\
\end{align*}
\ENDFOR
\end{algorithmic}
\end{algorithm}
Furthermore, Algorithm \ref{alg:batch_DISTA} can be  reformulated in an online fashion as illustrated in Algorithm \ref{alg:online_DISTA}. We denote this  online version as O-DISTA.
\begin{algorithm}  
\caption{O-DISTA for Problem \eqref{distributed_t}} \label{alg:online_DISTA}
\begin{algorithmic}[1]
\STATEx {\bf{input}}: $\lambda>0$, $\tau>0$; for each $v\in\V$,  $x_{v,0}=0$; at time $t$, $Q_{v,t}$,  $\phi_{v,t}$
\STATEx {\bf{output}}: in $[t,t+1)$, for each $v\in\V$, an estimate $x_{v,t+1}$ of $\xtrue_t$ 

\FOR{$t=0,\dots,T$}
\FOR{$h=0,\dots,r$}
\STATE For each $v\in\V$, $\ix_{v,0}=x_{v,t}$
\STATE If $h$ is even, for any $v\in\V$,
     \begin{align*}
c_{v,h+1}&=\overline{\ix}_{v,h}\\
 \ix_{v,h+1}&=\ix_{v,h}
\end{align*}
\STATE If $h$ is odd, for any $v\in\V$,
\begin{align*}
&c_{v,h+1}=c_{v,h}\\
& {\ix}_{v,h+1}=\soft_{\lambda\tau/2}\left[\frac{ {\ix}_{v,h}+\overline{c}_{v,h} - \tau  Q_{c,t} \ix_{v,h} -\tau\phi_{v,t}}{2}\right]
\end{align*}
\ENDFOR
\STATE For each $v\in\V$, $x_{v,t+1}=\ix_{v,r}$
\ENDFOR
\end{algorithmic}
\end{algorithm}
\section{Dynamic Regret Analysis for O-DISTA}\label{sec:dr_dista}
In this section, we analyze the dynamic regret for O-DISTA, by extending the results in \cite[Section IV]{fox18cdc} to the distributed case. We start by proving the contractivity.
\begin{lemma}\label{contr_d}
For each $t$, let  $\tau\leq\min_{v\in\V}\left\|A_{v,t}\right\|_2^{-2}$. If $X(t)\in\R^{n,\cardV}$ is the sequence produced by O-DISTA, and $\Xmin_t=(\xmin_t,\dots,\xmin_t)\in\R^{n,\cardV}$ is the minimizer of \eqref{distributed_t}, then
\begin{equation}
  \|X_{t+1}-\Xmin_t\|_F\leq \left(\frac{1+ \theta_{\tau}}{2}\right)^{r/2} \|X_{t}-\Xmin_t\|_F
 \end{equation}
 where $\|\cdot\|_F$ denotes the Frobenious norm, and $\theta_\tau=\max_{v,t}\|I-\tau Q_{v,t}\|_2^2<1$.
\end{lemma}
\begin{proof}
For an individual $v\in\V$, if $r=1$,
\begin{equation}
\begin{split}
&\left\|x_{v,t+1}-\xmin_t\right\|^2_2=\frac{1}{4}\left\|(I-\tau Q_{v,t})(x_{v,t}-\xmin_t)+ \overline{\overline{ x}}_{v,t}- \xmin_t  \right\|_2^2
\end{split}
 \end{equation}
 where  $\overline{\overline{x}}_{v,t}=\frac{1}{d^2}\sum_{w\in\Nv}\sum_{u\in\Nw}x_{u,t}$; notice that we use the fact that $\Xmin_t$ is a fixed point and that $\overline{\xmin_t}=\xmin_t$.
 
By applying the Cauchy-Schwarz inequality 
$\left(\sum_{i=1}^n a_i\right)^2\leq n\sum_{i=1}^n a_i^2$, we compute the following bound:
 \begin{equation}
 \begin{split}
  \|x_{v,t+1}-\xmin_t\|^2_2&\leq \frac{1}{2}\left\|I-\tau Q_{v,t}\right\|_2^2\left\|x_{v,t}-\xmin_t\right\|_2^2+\\
  &+\frac{1}{2}\frac{1}{d^2}\sum_{w\in\Nv}\sum_{u\in\Nw}\|x_{u,t}-\xmin_t\|_2^2.
\end{split}
 \end{equation}
 By summing over  $v\in\V$ and by exploiting the $d$-regularity of the graph, we obtain
 \begin{equation}
\|X_{v,t+1}-\Xmin_t\|_F^2\leq \frac{\theta_{\tau}+1}{2}\left\|X_{t}-\Xmin_t\right\|_F^2.
 \end{equation}
 The extension to $r>1$ is straightforward.
\end{proof}
By exploiting the contractivity, the following result can be proven.
Let $$\Delta_t:=\|\Xmin_t-\Xmin_{t-1}\|_F.$$
\begin{lemma}\label{lem:somme_d} For any $t=1,\dots,T$,
	\begin{equation*}
	\begin{split}
	(a)~~&\sum_{t=1}^T \left\|X_t -\Xmin_t\right\|_2\leq c_1 + c_2 \sum_{t=1}^T \Delta_t\\
	(b)~~&\sum_{t=1}^T \left\|X_t -\Xmin_t\right\|^2\leq c_3 + c_4 \sum_{t=2}^T \Delta_t^2+c_5 \sum_{t=1}^T \Delta_t\\
	\end{split}
	\end{equation*}
	with constants $c_i>0$, $i=1,\dots,5$.
\end{lemma}
The proof is a straightforward extension of the proof of \cite[Lemma 2]{fox18cdc}, and omitted for brevity.

\begin{lemma}\label{lem:surro_d} For each $t=1,\dots,T$,
\begin{equation}
F_{t-1}(X_{t})-F_{t-1}(\Xmin_{t-1})\leq \frac{2}{\tau}\left\| X_{t-1}-\Xmin_{t-1}\right\|_F^2. 
\end{equation}
\end{lemma}
\begin{proof}
 Let us consider the intermediate variables $\mathring{X}_h$, $h=0,\dots,r$, between $t-1$ and $t$, starting from $\mathring{X}_0=X_{t-1}$. 
 Since at each iteration $F_{t-1}$ is decreased by O-DISTA, and given the properties of the surrogate functional in \eqref{distributed_t_surro}, we have: 
\begin{equation*}
 \begin{split}
 &F_{t-1}(X_{t})  \leq F_{t-1}(\mathring{X}_{1}) = F_{t-1}(\mathring{X}_{1}, \overline{\mathring{X}}_{1}, \mathring{X}_{1})\\
 &\leq F_{t-1}(\mathring{X}_{1}, \overline{X}_{t-1}, X_{t-1}) \leq F_{t-1}(\Xmin_{t-1}, \overline{X}_{t-1}, X_{t-1}).
 \end{split}
 \end{equation*}
Therefore,
\begin{equation*}
 \begin{split}
 &F_{t-1}(X_{t}) - F_{t-1}(\Xmin_{t-1})\\
 &\leq F_{t-1}(\Xmin_{t-1}, \overline{X}_{t-1}, X_{t-1}) - F_{t-1}(\Xmin_{t-1})\\
 &\leq\frac{1}{\tau}\sum_{v\in\V}\left[\frac{1}{d} \sum_{w\in\mathcal{N}_v}\left\|\overline{x}_{w,t-1}-\xmin_{t-1} \right\|_2^2+ \left\|x_{v,t-1}-\xmin_{t-1} \right\|_2^2\right].
 \end{split}
 \end{equation*}
Then, the thesis is obtained by applying the Cauchy-Schwarz inequality and by exploiting the $d$-regularity of the graph.
\end{proof}
For any $X\in\R^{n,\cardV}$, let us define 
\begin{equation}
\begin{split}
D_{t}(X)&:=F_{t}(X)-F_{t-1}(X)\\
\Delta_{Q_{v,t}}&:=\left\| Q_{v,t}-Q_{v,t-1}\right\|_2\\
\Delta_{\phi_{v,t}}&:=\left\| \phi_{v,t}-\phi_{v,t-1}\right\|_2. 
\end{split}
\end{equation}
\begin{assumption}\label{ass:relaxed_boundedness}
We assume that $\sup_{v,t}\Delta_{Q_{v,t}}$ and $\sup_{v,t}\Delta_{\phi_{v,t}}$ are bounded. Moreover, if $\Delta_{Q_{v,t}}\neq 0$ for at least one $v$, we assume that $\sup_t \|\xmin_t\|$ is bounded.
\end{assumption}
We notice that this assumption is weaker than Assumption \ref{ass:limited}: $\phi_{v,t}$ does not need to be bounded; if $Q_{v,t}$ is constant in time, $\xmin_t$ does not need to be bounded. A possible application is discussed  Section \ref{sub:indoor}. Briefly, O-DISTA requires a weaker boundedness assumption than O-DR as for O-DISTA the sequence  $F_t(X_t)$ is monotone decreasing,  and in particular Lemma \ref{lem:surro_d} holds, while this is not generally guaranteed for O-DR. 
\begin{lemma}\label{lem:D_t} For any $X\in\R^{n,\cardV}$ and $t=1,\dots,T$,
		\begin{equation}\label{unifcont2}
	D_t(X)-D_t(\Xmin_t)\leq \gamma_1 \|X-\Xmin_t\|_F +\gamma_2 \|X-\Xmin_t\|_F^2
	\end{equation} 
	where  $\gamma_1=\sqrt{\cardV} \sup_{v,t}(\Delta_{\phi_{v,t}}+\Delta_{Q_{v,t}}\|\xmin_t\|_2)$ and  $\gamma_2= \frac{1}{2}\sup_{v,t}\Delta_{Q_{v,t}}$.
\end{lemma}
\begin{proof}
Since 
\begin{equation*}
D_t(X)=\hspace{-0.1cm}\sum_{v\in\V}\left[\frac{1}{2}x_v^T (Q_{v,t}-Q_{v,t-1})x_v + (\phi_{v,t}-\phi_{v,t-1})^Tx_v\right] 
\end{equation*} 
 and by using \eqref{algebra}, we have:

\begin{equation*}
\begin{split}
&D_t(X)-D_t(\Xmin_t)=\\&
=\sum_{v\in\V}\left[\frac{1}{2}(x_v-\xmin_t)^T (Q_{v,t}-Q_{v,t-1})(x_v+\xmin_t) \right]\\
&+ \sum_{v\in\V}(\phi_{v,t}-\phi_{v,t-1})^T (x_v-\xmin_t) \\
&\leq\hspace{-0.1cm}\sum_{v\in\V}\left[\frac{1}{2}\left\|x_v-\xmin_t\right\|_2 \Delta_{Q_{v,t}}\left\|x_v+\xmin_t \right\|_2+\Delta_{\phi_{v,t}} \left\|x_v-\xmin_t\right\|_2\right]\\
&\leq \sum_{v\in\V}\frac{1}{2}\left\|x_v-\xmin_t\right\|_2^2 \Delta_{Q_{v,t}}+\\
& +\sum_{v\in\V}\left(\Delta_{\phi_{v,t}}+\Delta_{Q_{v,t}}\|\xmin_t\|_2 \right)\left\|x_v-\xmin_t\right\|_2.
\end{split}
\end{equation*}
%
Then, the thesis is obtained by applying Cauchy-Schwarz.
\end{proof}
Given these intermediate lemmas, we can now evaluate the dynamic regret.
\begin{theorem}\label{theo:regret_d}
The dynamic regret for O-DISTA  (Algorithm \ref{alg:online_DISTA}) has the following bound:
\begin{equation*}
\begin{split}
&\reg \leq \alpha_0+\alpha_1 \sum_{t=1}^T \Delta_t + \alpha_2 \sum_{t=1}^T \Delta_t
\end{split}
\end{equation*}
where $\alpha_i>0$, $i=0,1,2$.
\end{theorem}

\begin{proof}
We consider  the difference of losses $F_t(X_t)-F_t(\Xmin_t)$ and we add and subtract $F_{t-1}(X_t)$ to it. By using that $F_{t-1}(\Xmin_t)\geq F_{t-1}(\Xmin_{t-1})$, we obtain the following bound:
\begin{equation}\label{needname}
\begin{split}
&F_t(X_t)-F_t(\Xmin_t)\leq \\
&\leq F_t(X_t)-F_t(\Xmin_t)\pm F_{t-1}(X_t)+ F_{t-1}(\Xmin_{t})-F_{t-1}(\Xmin_{t-1})\\
&= D_t(X_t)-D_t(\Xmin_t)+F_{t-1}(X_t)-F_{t-1}(\Xmin_{t-1}).
\end{split}
\end{equation}
Then, by applying Lemma \ref{lem:surro_d} and Lemma \ref{lem:D_t}, the last expression is upper bounded by
\begin{equation}
\gamma_1 \|X_t-\Xmin_t\|_F +\gamma_2 \|X_t-\Xmin_t\|_F^2+\frac{2}{\tau}\left\| X_{t-1}-\Xmin_{t-1}\right\|_F^2. 
\end{equation}
The thesis is obtained by summing over $t=1,\dots,T$, and by applying Lemma \ref{lem:somme_d}.
\end{proof}

%
\section{Numerical results}\label{sec:nr}
 \begin{figure*}[h!]
	\centering
	\includegraphics[width=0.49\columnwidth]{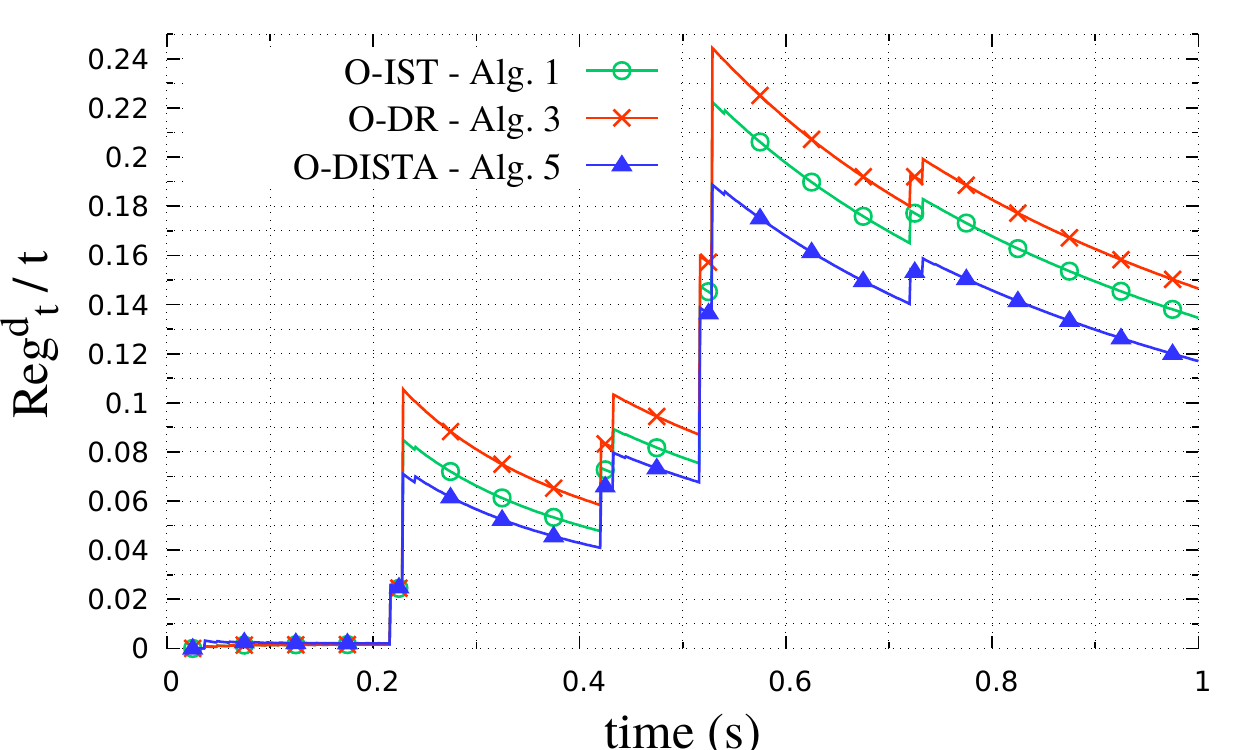}
	\includegraphics[width=0.49\columnwidth]{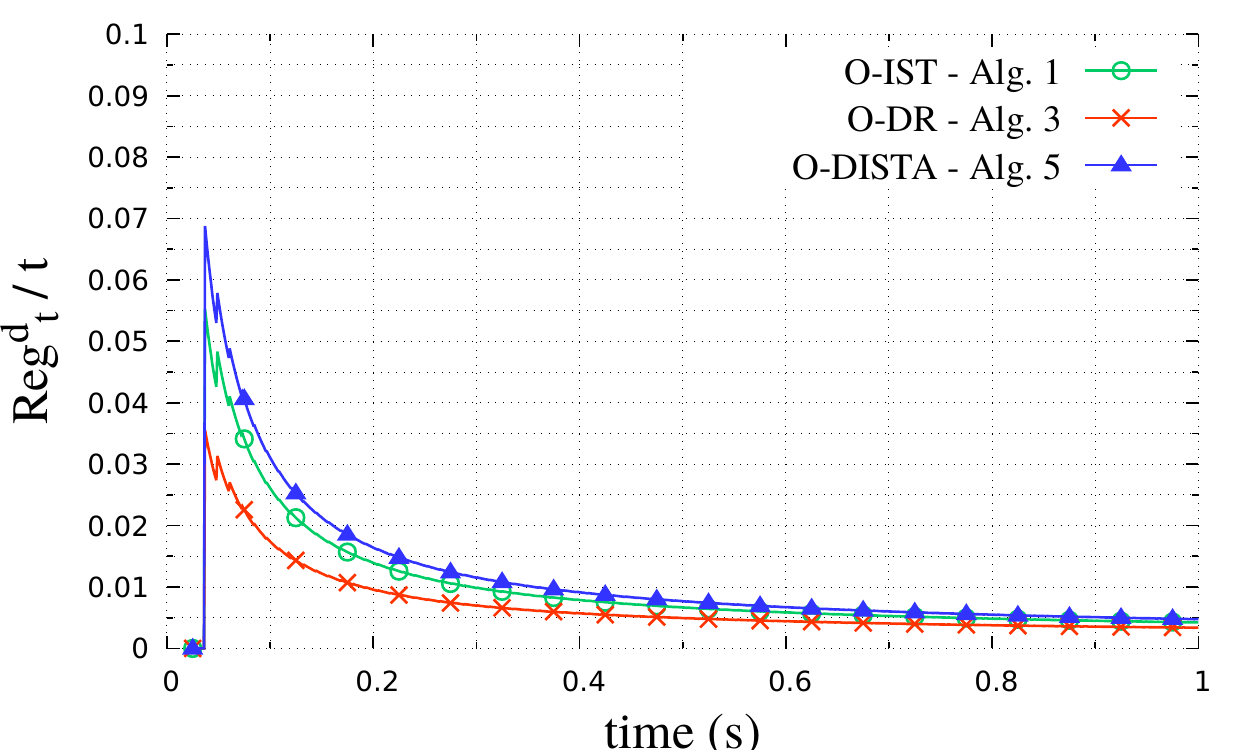}
	\includegraphics[width=0.49\columnwidth]{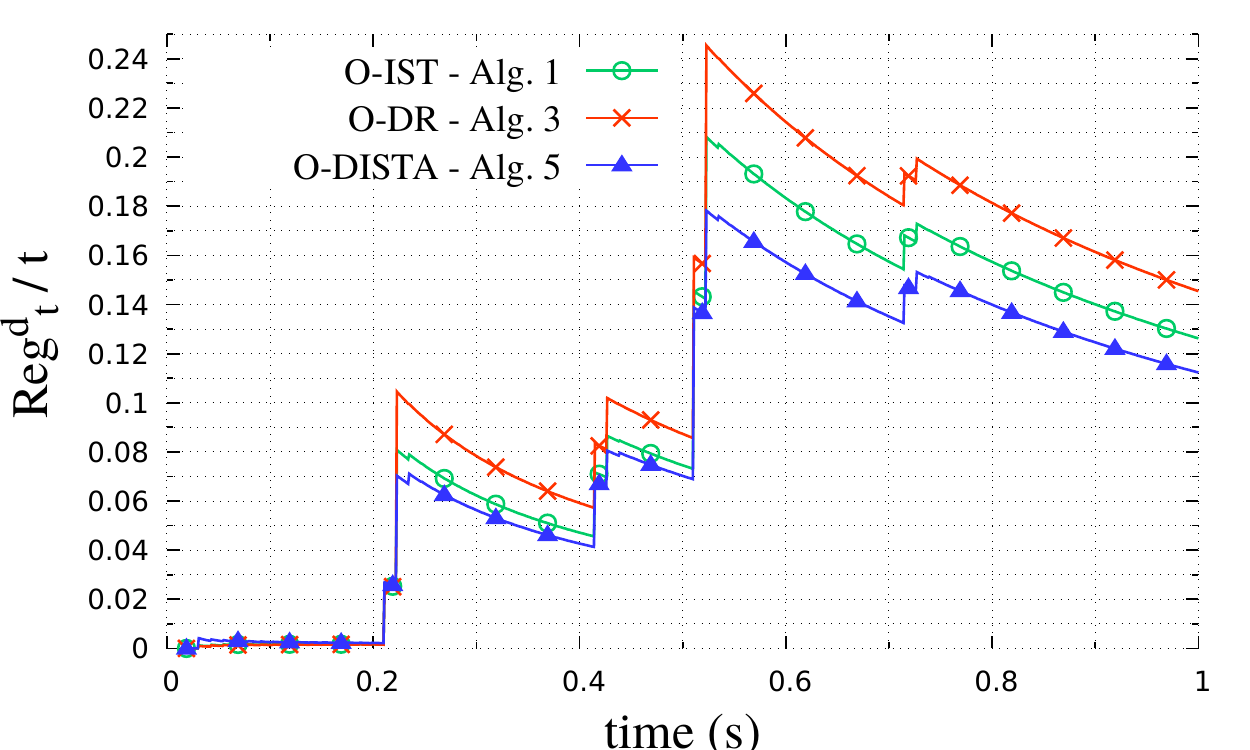}
	\includegraphics[width=0.49\columnwidth]{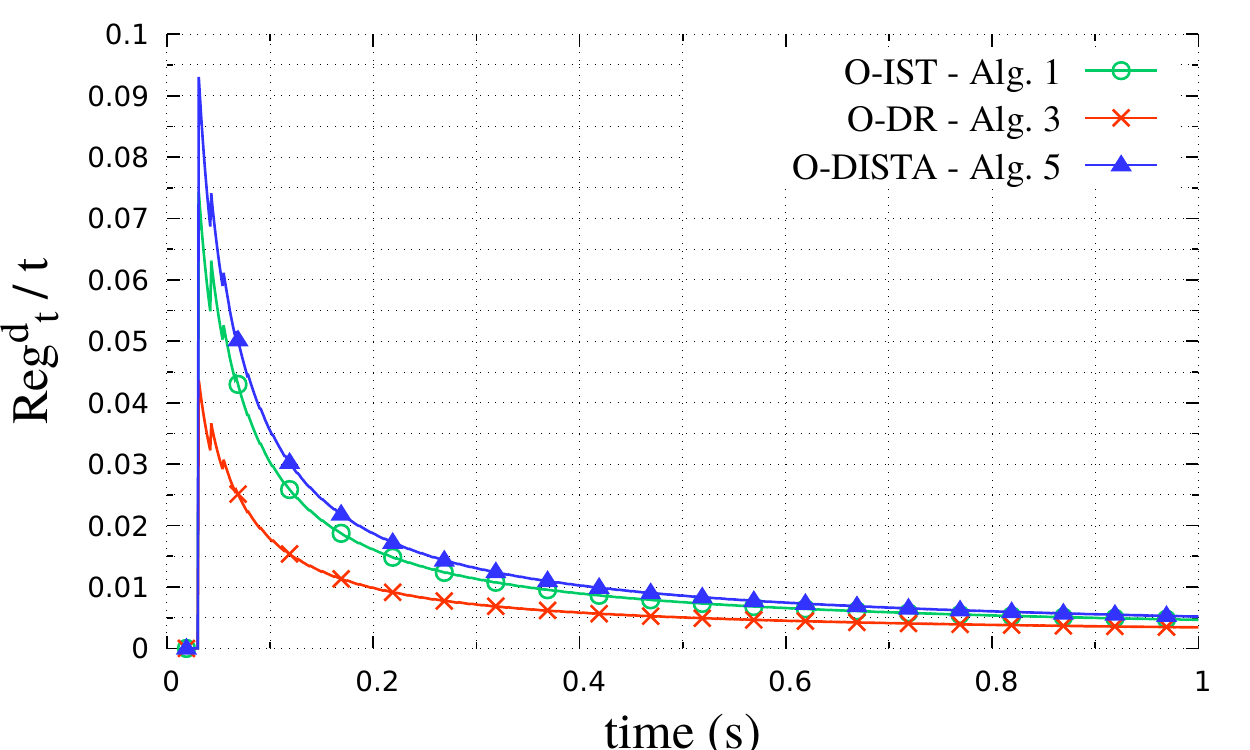}
	\caption{Dynamic regret. Left column: Experiment 1; right column: Experiment 2. From top to bottom: (1) $t_r=12$ ms, SNR=25dB; (2) $t_r=6$ ms, SNR=25dB.}
	\label{fig:s0}
\end{figure*}
In this section, we numerically analyze the algorithms O-IST, O-DR, and O-DISTA in two time-varying Elastic-net experiments.{\footnote{The code to reproduce the proposed experiments is available at https://github.com/sophie27/Sparse-Time-Varying-Optimization.} The first problem is an instance of online identification of time-varying linear systems; the second problem is a practical example of moving target tracking based on the received signal strength (RSS), which has applications, e.g., in indoor monitoring and surveillance.
\subsection{Online compressed system identification}\label{sub:csi}
Online compressed system identification refers to the online estimation of the parameters of a time-varying system from compressed measurements. Specifically, we consider a time-varying autoregressive model with an exogenous input (TVARX), whose input-output relationship is as follows: $y_t=\sum_{p=1}^P a_{p,t}y_{t-p}+\sum_{q=1}^Q b_{q,t}u_{t-q}+e_t$, 
 where $u_t,y_t\in\R$ respectively are the measurable input and output; $e_t\in\R$ is the measurement error; $a_{p,t}, b_{q,t}\in\R$ are the time-varying parameters to be estimated. 

 The dimensions $P$ and $Q$ are assumed to be unknown, therefore we initially set sufficiently large bounds $\widehat{P}$ and $\widehat{Q}$  for them and then we look for a parsimonious model using the $\ell_1$-norm to promote sparsity.

As in \cite{fox18cdc}, we iteratively collect groups of $m$ measurements $\mathbf{y}_t:=(y_t,\dots, y_{t+m})^T$. The measurements are compressed, that is, we choose $m<\widehat{P}+\widehat{Q}$; a  Gaussian measurement noise with SNR$=$25dB is added.

It is easy to check that we can define $A_t\in\R^{m, \widehat{P}+\widehat{Q}}$ as follows:
\begingroup\makeatletter\def\f@size{9}\check@mathfonts
\begin{equation*}
\left(\begin{array}{cccccc}
		y_{t-1}&\cdots&y_{t-P}&u_{t-1}&\cdots&u_{t-Q}\\
		y_{t}&\cdots&y_{t-P+1}&u_{t}&\cdots&u_{t-Q+1}\\
		\vdots&&&&&\vdots\\
		y_{t+m-1}&\cdots&y_{t+m-P}&u_{t+m-1}&\cdots&u_{t+m-Q}\\
	\end{array}\right).
\end{equation*}
\endgroup
We revisit the TVARX(1,1) example considered in \cite{li11,fox18cdc}. The input-output equation is $y_t=a_{1,t}y_{t-1}+b_{1,t}u_{t-1}+e_t$, with $P=Q=1$. Assuming $P$ and $Q$ unknown, we initially overestimate them as $\widehat{P}=\widehat{Q}=10$. Thus, our goal is to  track a time-varying sparse vector $\xtrue_t=(a_{1,t},\dots,a_{10,t},b_{1,t},\dots, b_{10,t})\in\R^{n}$, $n=20$, with sparsity $k=2$ and constant support, given linear observations, as in \eqref{acquisition}. The Elastic-net model \eqref{eq:elasticnet} is efficient to tackle this problem, as shown in \cite{fox18cdc}. By cross-validation, we set $\lambda=10^{-2}$ and $\mu=10^{-6}$. We remark that, for stability purpose, it makes sense to assume $A_t$ and $\xtrue_t$ bounded, which matches with Assumption \ref{ass:limited}.
\begin{figure*}[h!]
	\centering
	\includegraphics[width=0.49\columnwidth]{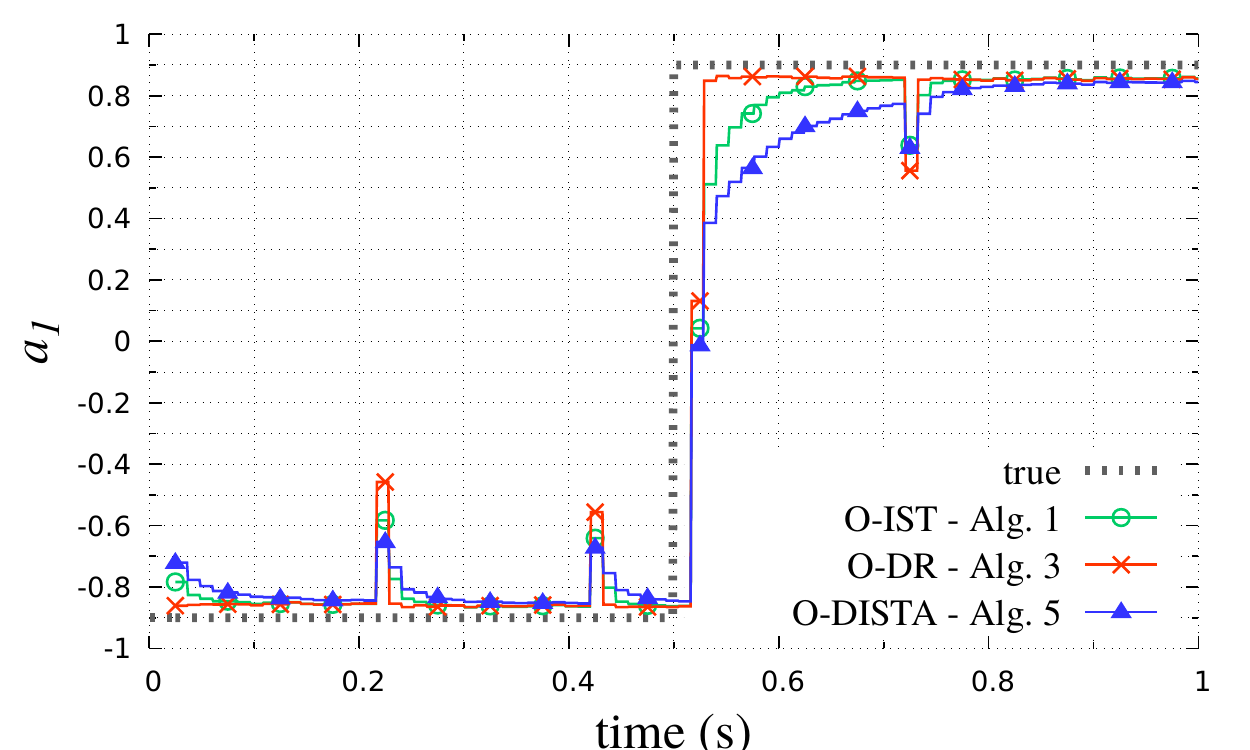}
	\includegraphics[width=0.49\columnwidth]{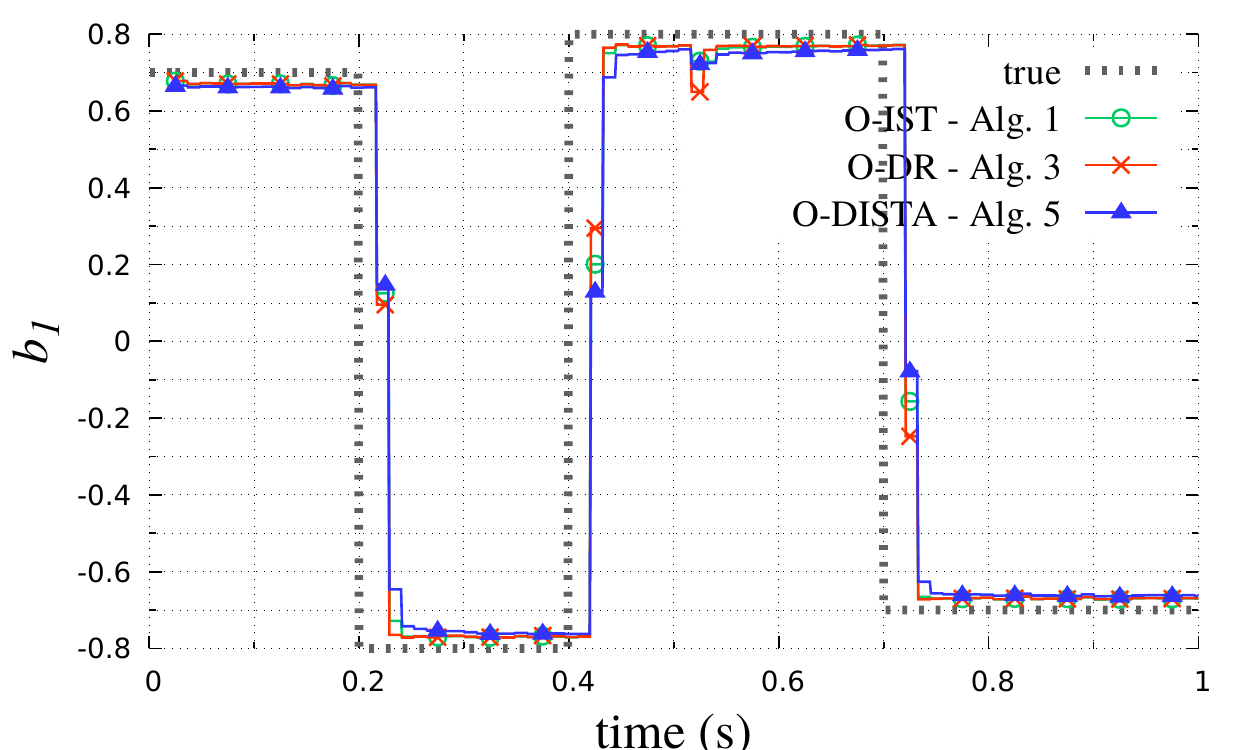}
	\includegraphics[width=0.49\columnwidth]{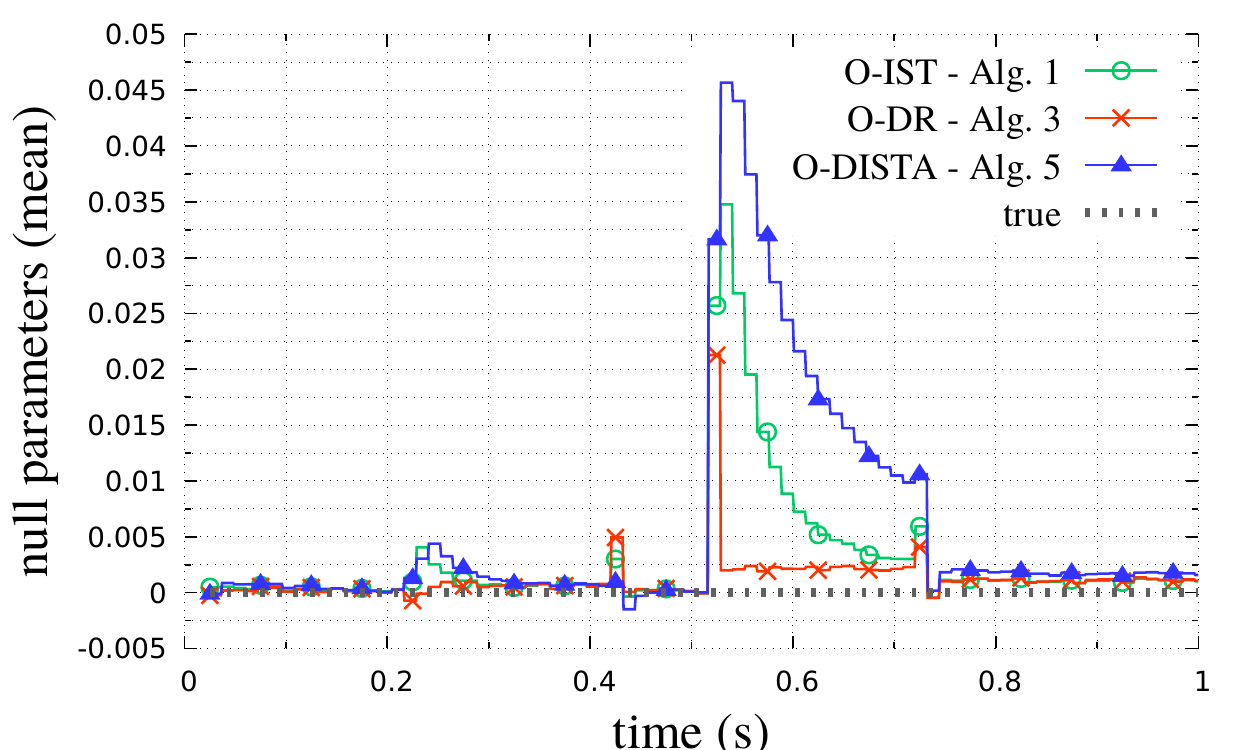}
	\includegraphics[width=0.49\columnwidth]{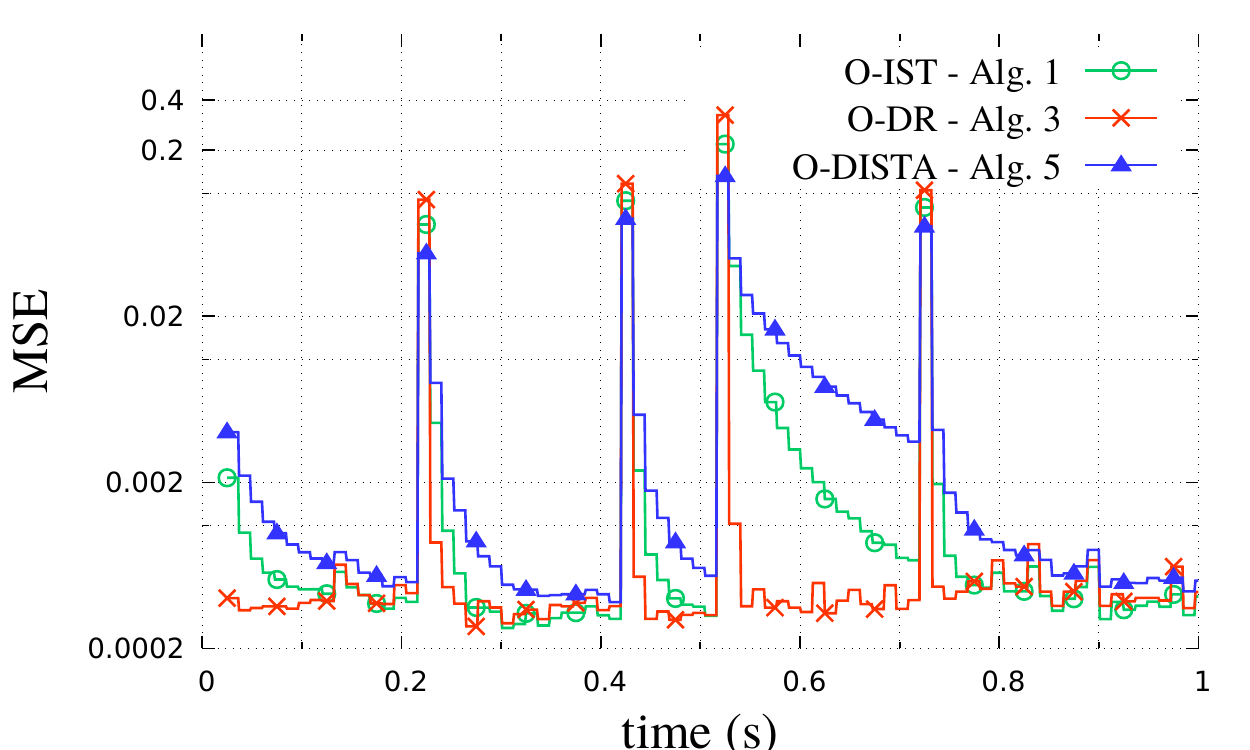}
	\caption{Experiment 1, SNR=25dB, $t_r=12$ ms. From left to right, averaged estimates of $a_{1,t}$, $b_{1,t}$, null parameters; mean square error.}
	\label{fig:s1}
\end{figure*}
\begin{figure*}[h!]
	\centering
	\includegraphics[width=0.49\columnwidth]{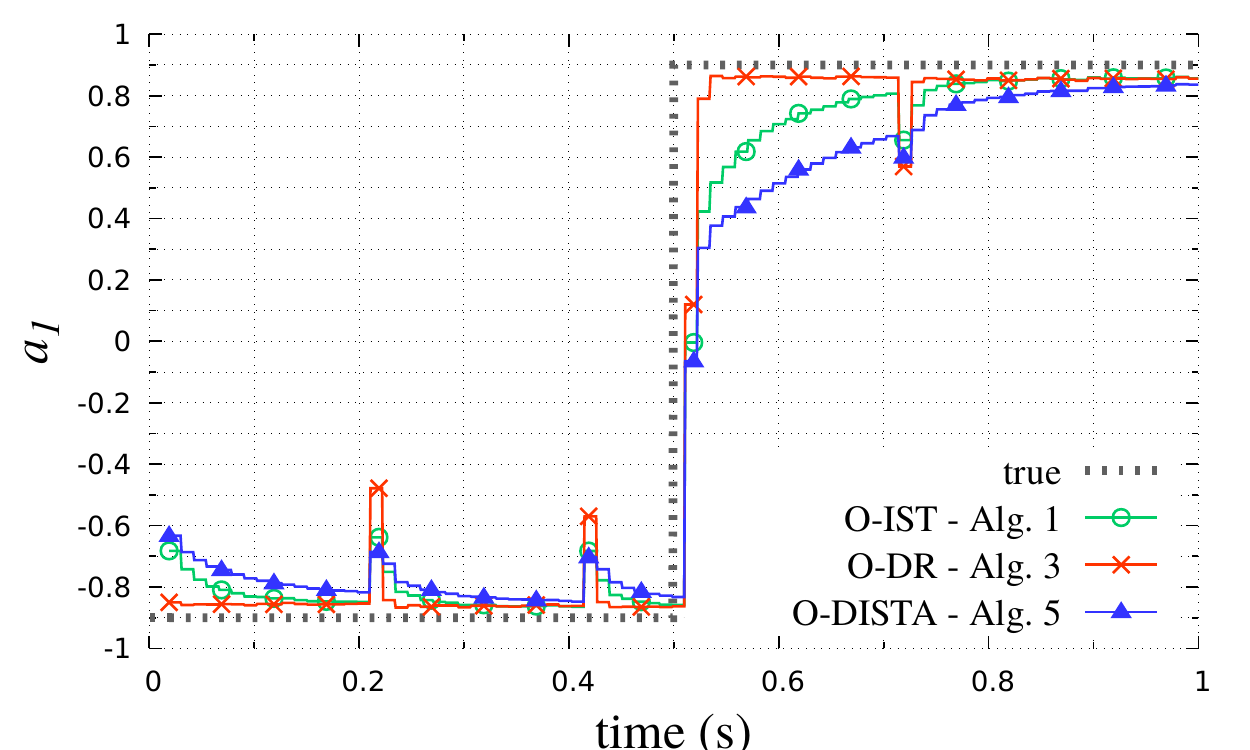}
	\includegraphics[width=0.49\columnwidth]{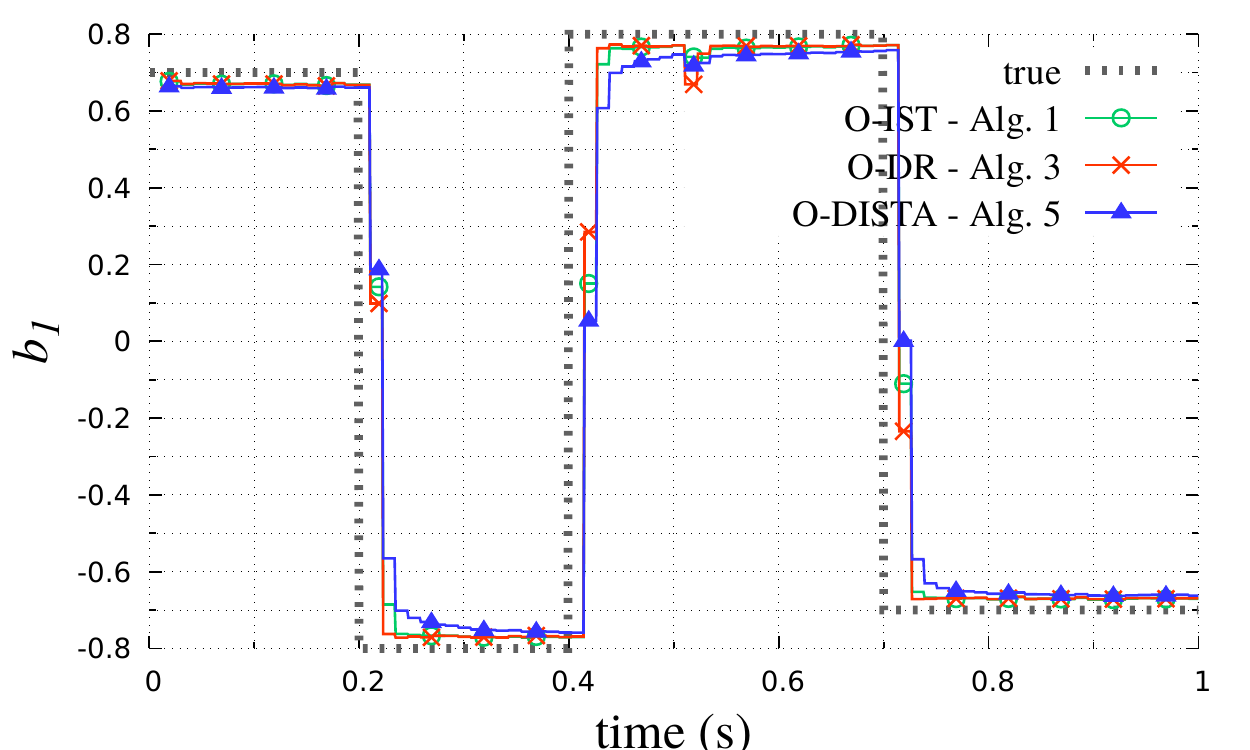}
	\includegraphics[width=0.49\columnwidth]{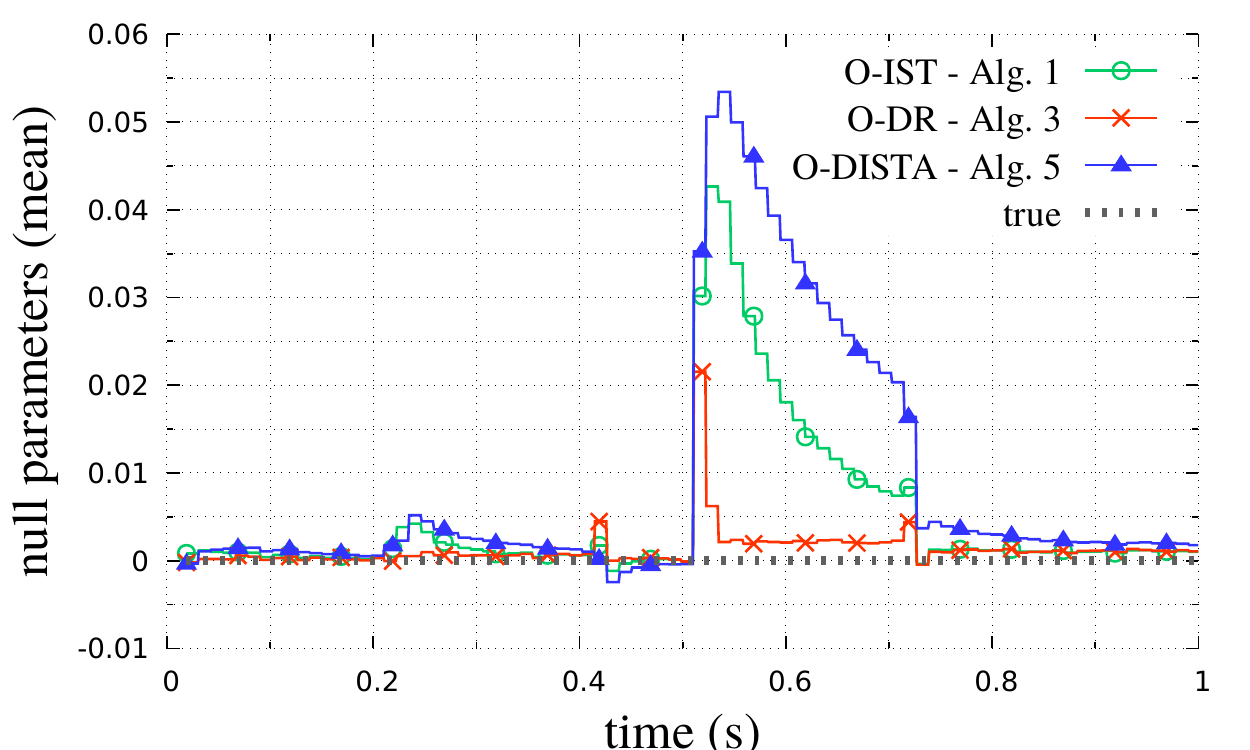}
	\includegraphics[width=0.49\columnwidth]{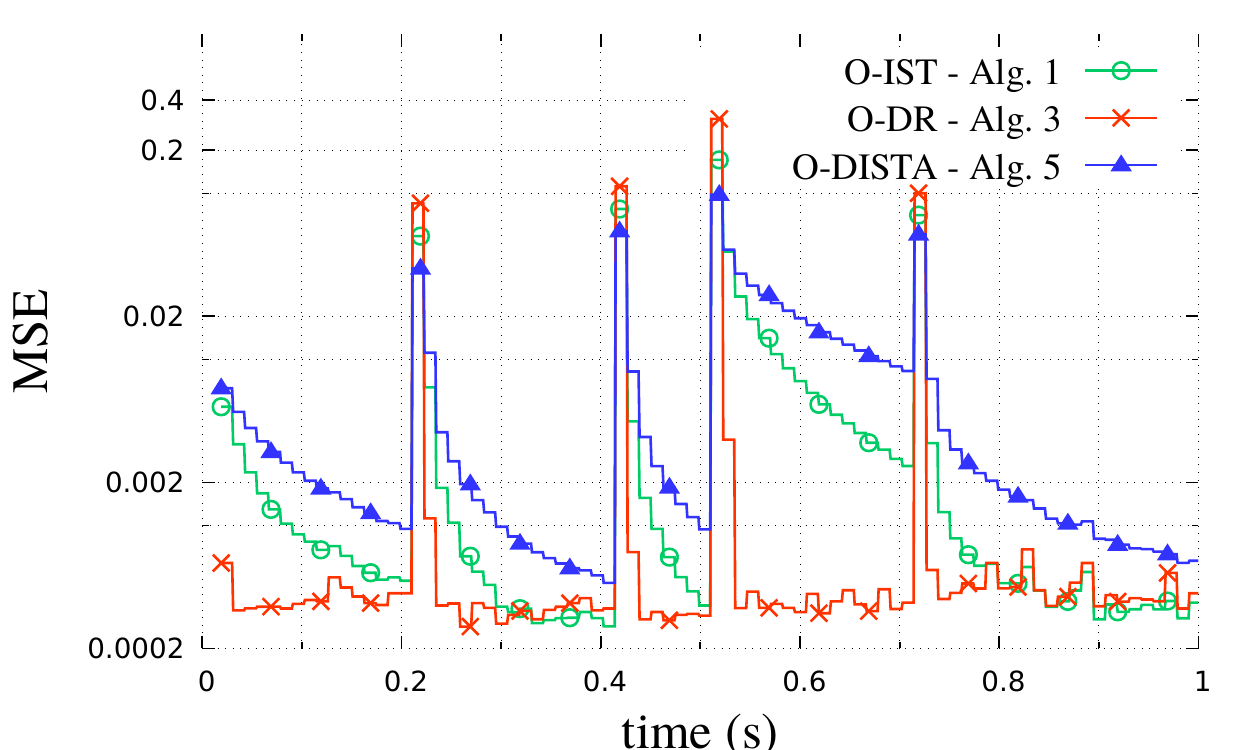}
	\caption{Experiment 1, SNR=25dB, $t_r=6$ ms. From left to right, averaged estimates of $a_{1,t}$, $b_{1,t}$, null parameters; mean square error.}
	\label{fig:s2}
\end{figure*}
\begin{figure*}[h!]
	\centering
	\includegraphics[width=0.49\columnwidth]{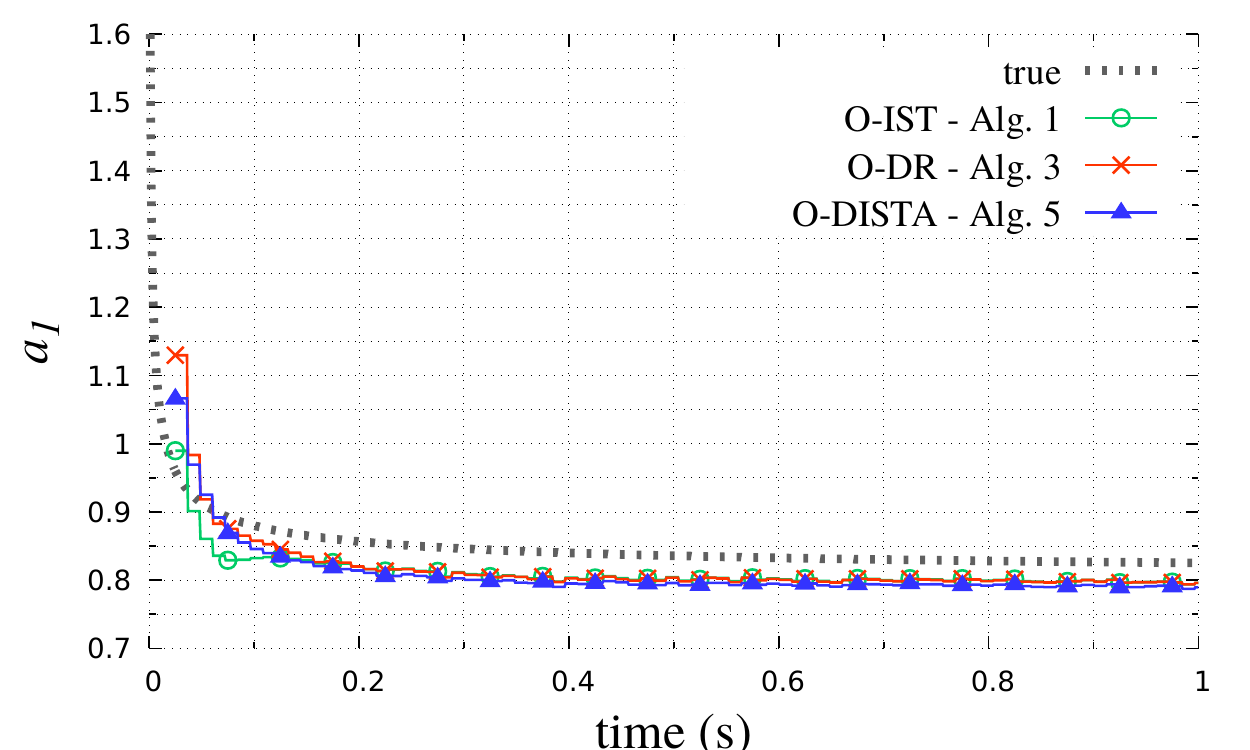}
	\includegraphics[width=0.49\columnwidth]{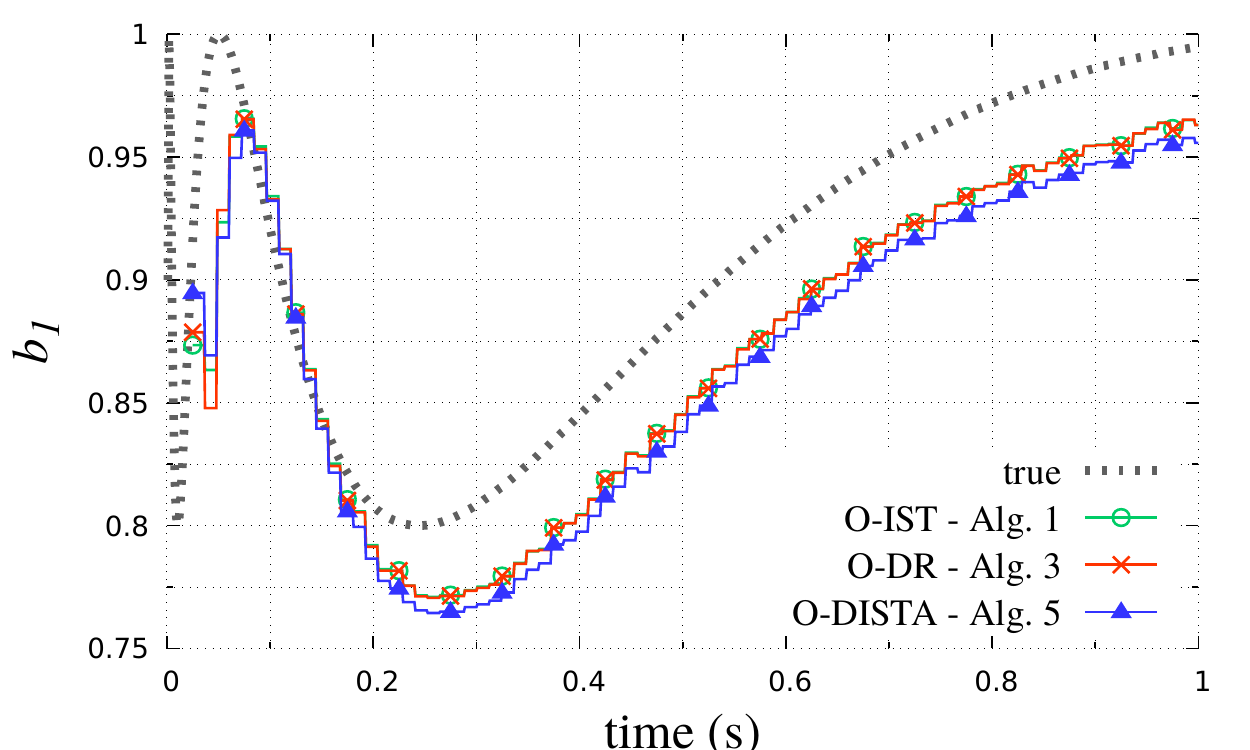}
	\includegraphics[width=0.49\columnwidth]{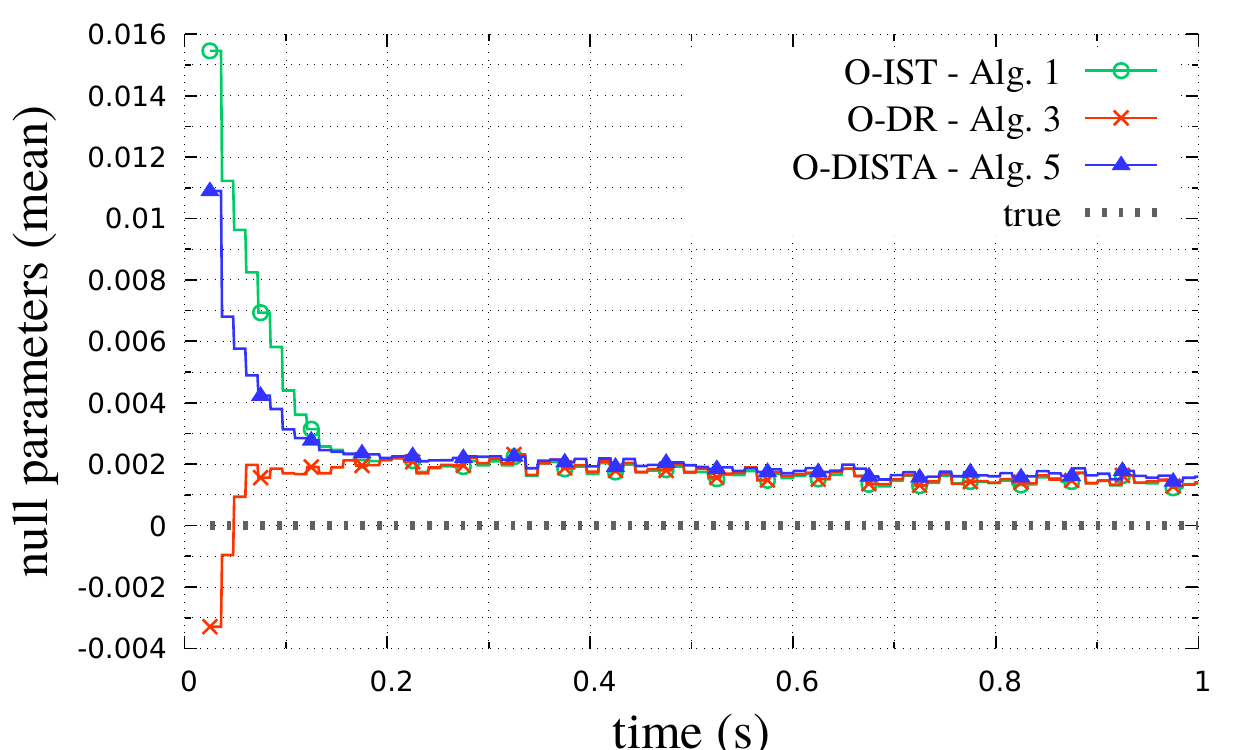}
	\includegraphics[width=0.49\columnwidth]{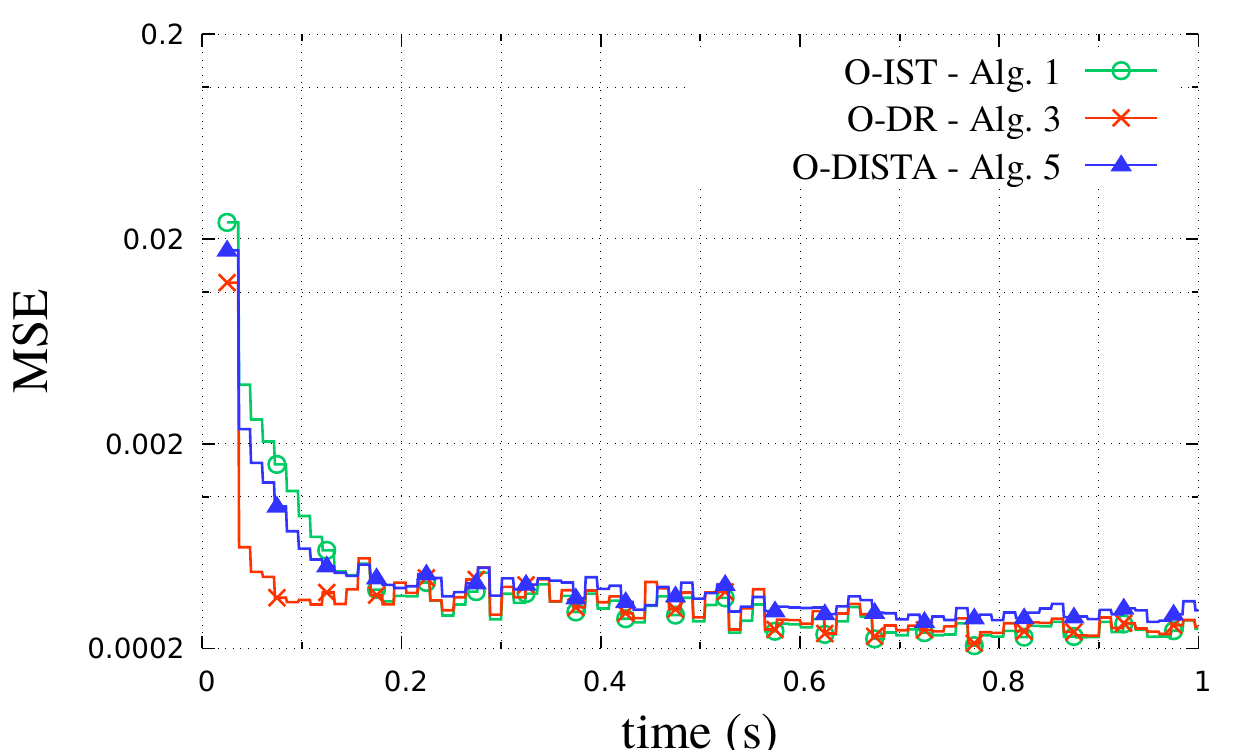}
	\caption{Experiment 2, SNR=25dB, $t_r=12$ ms. From left to right, averaged estimates of $a_{1,t}$, $b_{1,t}$, null parameters; mean square error.}
	\label{fig:s3}
\end{figure*}
\begin{figure*}[h!]
	\centering
	\includegraphics[width=0.49\columnwidth]{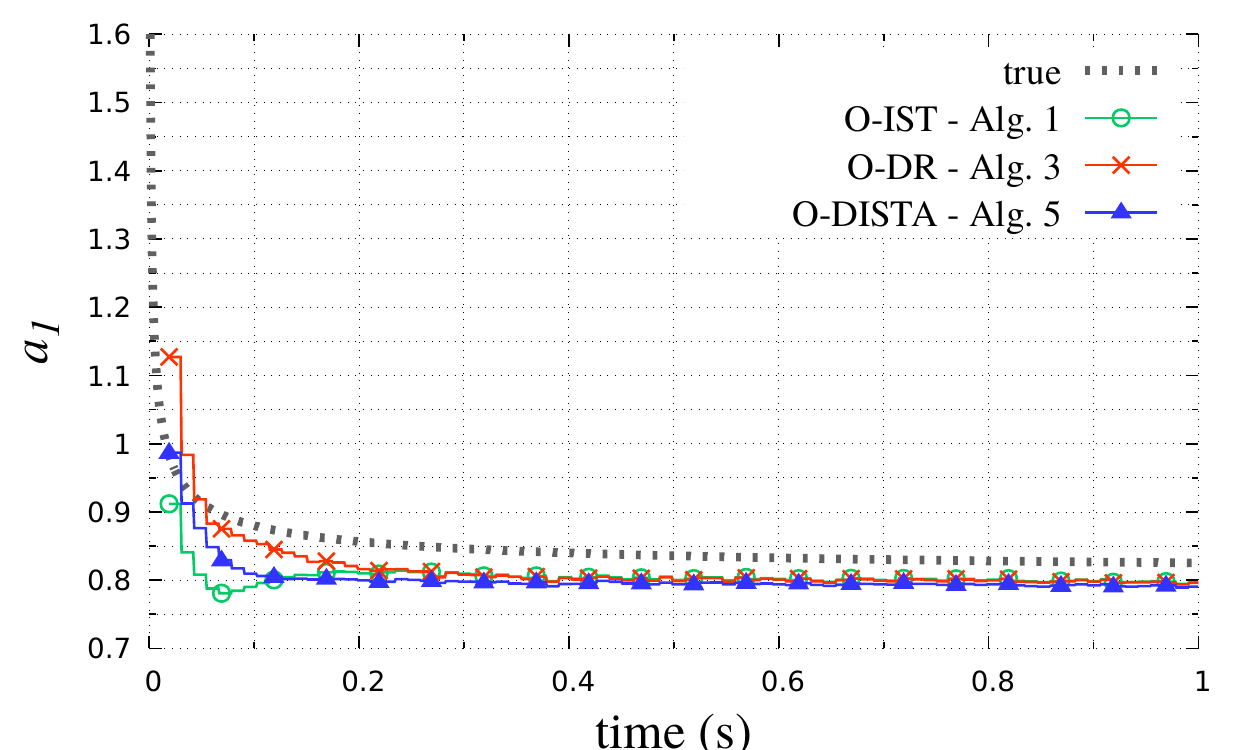}
	\includegraphics[width=0.49\columnwidth]{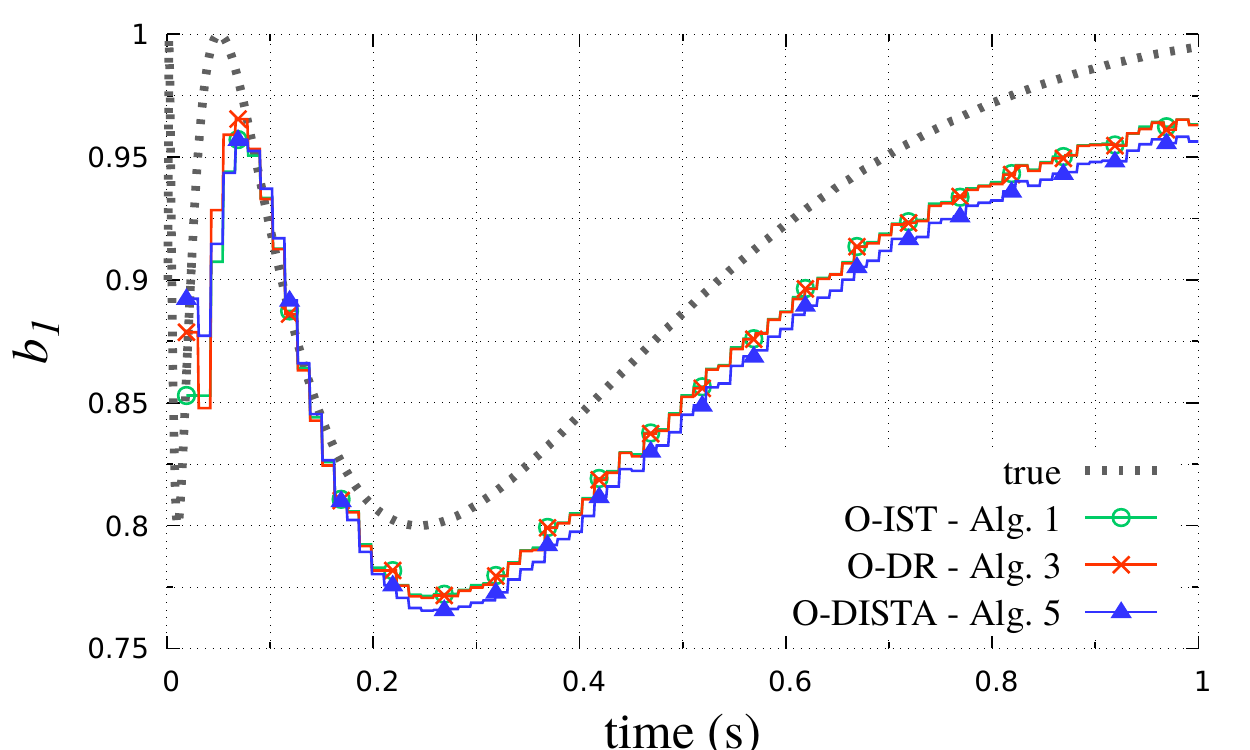}
	\includegraphics[width=0.49\columnwidth]{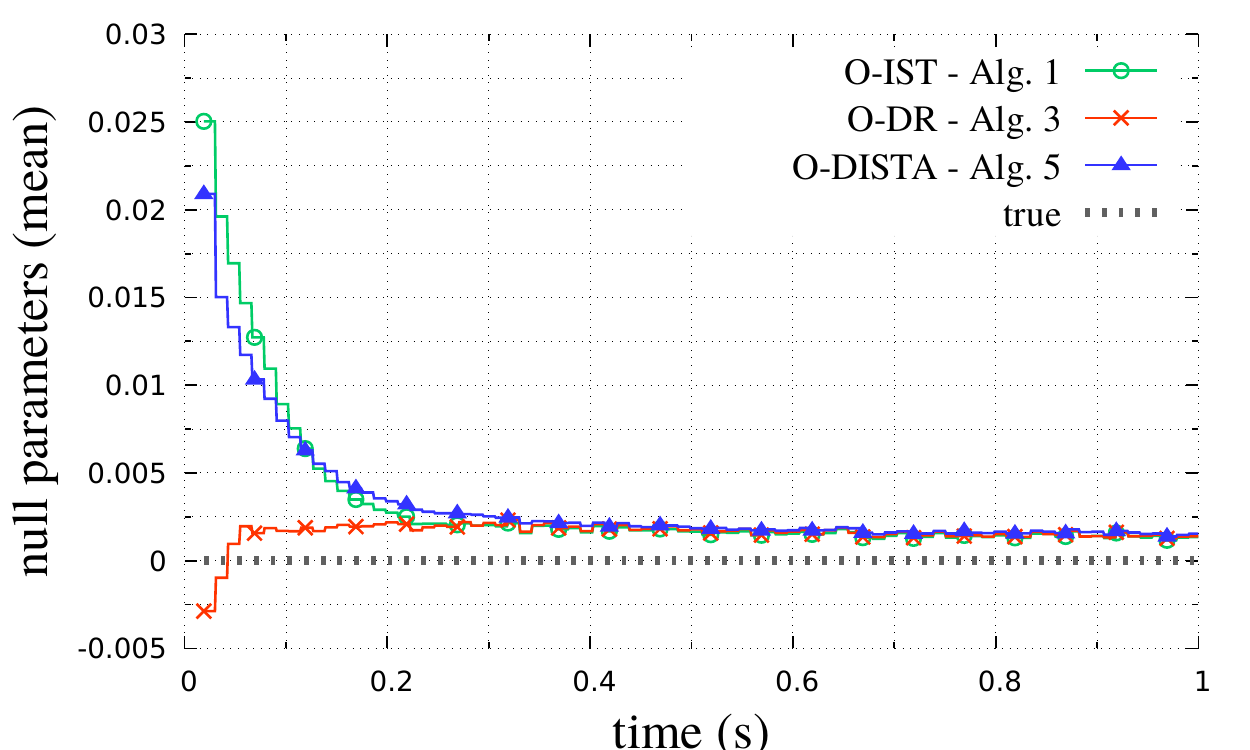}
	\includegraphics[width=0.49\columnwidth]{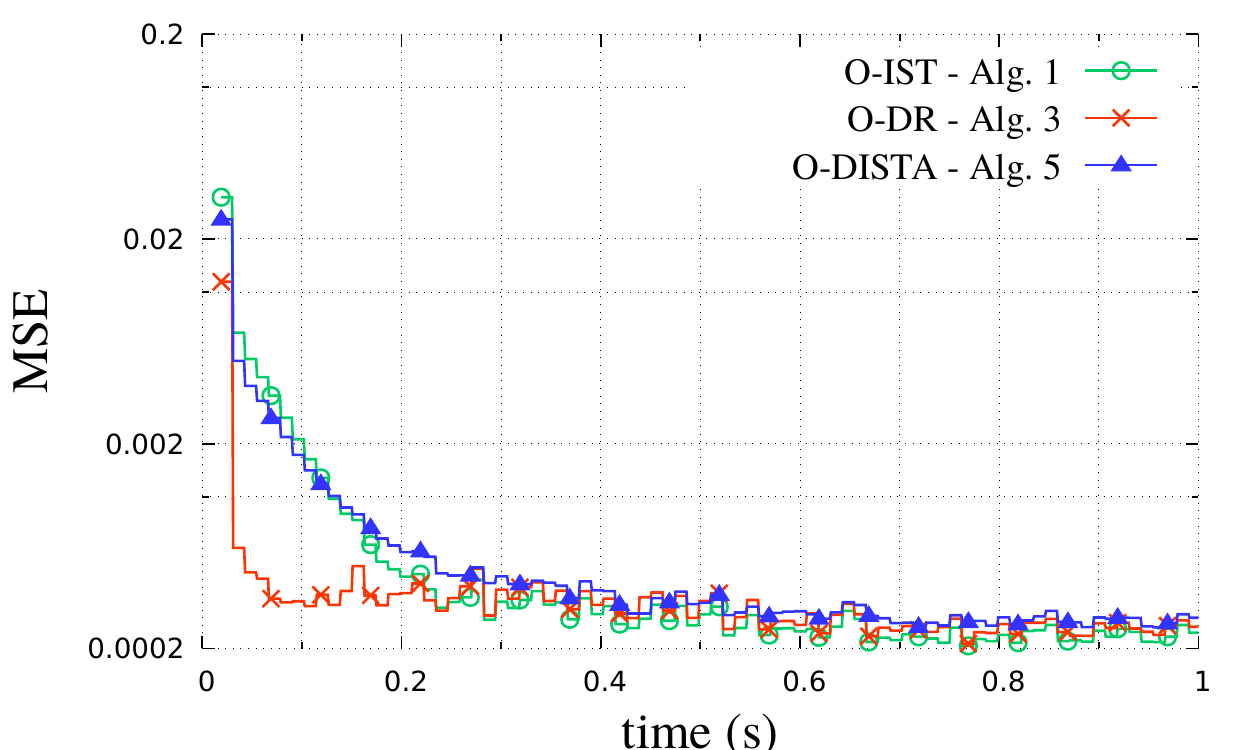}
	\caption{Experiment 2, SNR=25dB, $t_r=6$ ms. From left to right, averaged estimates of $a_{1,t}$, $b_{1,t}$, null parameters; mean square error.}
	\label{fig:s4}
\end{figure*}

Information on sparsity and support is assumed to be unknown and not exploited in the estimation. We consider a time horizon of $1$ second and sampling frequency of 1000 Hz. Two experiments are conducted. In the first one, $a_{1,t}$ and $b_{1,t}$ are step-wise constant with few abrupt changes, see \cite{li11}. Specifically, we set:

\textbf{Experiment 1:}

\begingroup\makeatletter\def\f@size{9}\check@mathfonts
\begin{equation}
a_1(t)=\left\{\begin{aligned}
&-0.9 \text { if } t<0.5\\
&0.9 \text{ otherwise;}
\end{aligned}\right.
~~b_1(t)=\left\{\begin{aligned}
&0.7 \text { if } t<0.2\\
&-0.8 \text { if } 0.2 \leq t<0.4\\
&0.8 \text { if } 0.4 \leq t<0.7\\
&-0.7 \text{ otherwise.}\\
\end{aligned}\right.
\end{equation}
\endgroup
In the second experiment, instead, we test a case of smoothly time-varying parameters.

\textbf{Experiment 2:}
\begin{equation*}
a_1(t)=0.8\left(1+\frac{1}{\sqrt{t}}\right);~~b_1(t)=0.9+0.1 \sin \left(2 \log t\right).
\end{equation*}
These parameters are chosen so that $a_1(t)$ is decreasing to zero, with convergent path length $\sum_{t=1}^{T} |a_1(t+1)-a_1(t)|$, while $b_1(t)$ is oscillating, with sublinear path length $\sum_{t=1}^{T} |b_1(t+1)-b_1(t)|$, of order $\log T$. We notice that the path length was defined above on the optimal points $\xmin_t$, while here we are evaluating it on $\xtrue_t$; however, $\xmin_t$ is expected to be a good approximation of $\xtrue_t$, then the two path lengths are somehow equivalent, see \cite[Corollary 1]{fox18cdc}.
The input components are drawn from a standard Gaussian distribution, and are  periodic with period $m$. 
We set $m=12$, which corresponds to a rate compression $\frac{m}{n}=\frac{3}{5}$. This implies a delay of 12 ms to acquire the of set measurements plus the run time. To prevent an accumulation of delay, the run time of the algorithm must not exceed 12 ms, so that the algorithm processes the acquired measurements while the successive measurements are being acquired. Then, we  iterate the algorithms until a prefixed maximum run time $t_r\leq m$ ms (notice that this approach is different from that of \cite{fox18cdc}, where a maximum number of iterations was set). In our simulations, we test $t_r=12$ ms and $t_r=6$ ms.  

For O-IST and O-DISTA, sufficient conditions on the parameter $\tau$ have been theoretically provided; specifically, for each $t$, $\tau\leq \left\|A_{t}\right\|_2^{-2}$ for O-IST (see Assumption 1 in \cite{fox18cdc}), and  $\tau\leq\min_{v\in\V}\left\|A_{v,t}\right\|_2^{-2}$ for O-DISTA, see Lemma \ref{contr_d}. In these experiments, we assume to ignore these lower bounds, and we set $\tau$ at each time step: specifically, we use  $\tau=2\left\|A_t\right\|_2^{-2}$ for O-IST, while for O-DISTA, each node computes its own $\tau$ as $2\left\|A_{v,t}\right\|_2^{-2}$. These values are observed to keep the convergence properties of the algorithms in practice. For O-DISTA, we consider a $3$-regular ring topology: there are 4 nodes, each of them taking 3 measurements and communicating with 2 neighbors.

The results shown in figures \ref{fig:s0}-\ref{fig:s4} are averaged over 200 random runs. 
In Figure \ref{fig:s0}, we show the evolution of the mean dynamic regret $\mathrm{\mathbf{Reg}}^d_t/t$. Based on our theoretical results, we expect that $\mathrm{\mathbf{Reg}}^d_t/t$ decreases to zero (that is, the dynamic regret is sublinear), when the system is static or when it evolves with sublinear path length. This behavior is confirmed by numerical simulations and can be appreciated  in  Figure \ref{fig:s0}.

In figures  \ref{fig:s1}-\ref{fig:s4}, we illustrate more details for each experiment. Specifically, we provide four graphs, respectively  depicting the tracking of $a_1(t)$,  $b_1(t)$, and  null parameters, and the mean square error, defined as  MSE=$\frac{1}{P+Q}\sum_{s=1}^{T/m}\left\|\xtrue_{sm}-\widehat{x}_{sm}\right\|_2^2$, where $\xtrue_t=(a_{1,t},\dots, a_{P,t}, b_{1,t},\dots,b_{Q,t})^T$ and $\widehat{x}_t$ is the estimate (the mean estimate depicted for O-DISTA).

Concerning Experiment 1 (figures \ref{fig:s1}-\ref{fig:s2}), in general, all the implemented algorithms are able to track the true parameters. When the parameters jump between different values (this occurs at time instants 0.2, 0.4, 0.5, 0.7), all the estimates are affected by a sudden perturbation. O-DR is observed to  converge faster than O-IST and O-DISTA when the parameters are constant. After jumps, O-DR adapts faster to the new parameters. On the other hand, O-DR is locally more sensitive to jumps: its peaks in correspondence of jumps are more marked. As expected, the distributed nature of O-DISTA makes it a bit less prompt than the centralized algorithms.

 Moving from $t_r=12$ ms to $t_r=6$ ms, we obtain a slight worsening for all the algorithms; on the other hand, the response delay after jumps is reduced. O-DR performance is almost equal for $t_r=12$ ms and $t_r=6$ ms, which suggests that O-DR almost achieves convergence to the optimal point in 6 ms. 

Concerning Experiment 2 (figures \ref{fig:s3}-\ref{fig:s4}), similar considerations can be drawn. In addition, we observe that O-DR and O-IST have similar performance when the parameters are slowly varying (namely, for $t\geq200$ ms), while O-DR is more precise when the path length is higher ($t<200$ ms).
\subsection{Moving target tracking}\label{sub:indoor}
In the last decade, indoor localization of moving objects has been gaining attention for purposes such as monitoring and surveillance, tracking of products in manufacturing industrial lines, control of unmanned vehicles, and location-based services. While outdoor tracking is mature, due to satellites technologies, indoor tracking is still challenging, and a variety of methodologies are proposed for it, see \cite{dar15} for a complete overview.

A possible approach to indoor tracking is based on the distance estimation via RSS, which can be implemented in 
low-cost systems such as wireless sensor networks. RSS-positioning is often associated with CS techniques to obtain an accurate localization from few measurements, see \cite{fen09,fen12}.

In this experiment, we consider the CS model proposed in \cite{fen09} and we extend it to the dynamic case, by using an Elastic-net model. 
Specifically, we aim to track a moving target in a $25\times 25~\text{m}^2$ indoor area. The area is assumed to subdivided into square cells of side $1$ m; the target is well localized when the cell where it lies is identified (actually, the target is localized in the center of the cell). This is sparse problem since only one cell over $n=625$ is occupied at each time step. The distance is measured with $36$ sensor nodes, deployed according to a regular grid over the area, represented by the yellow points in Figure \ref{fig:s5}.  Measurements are linearly obtained as $y_t=A\xtrue_t+$ noise, $y_t\in\R^m$,  through a dictionary $A\in\R^{m,n}$ built in a training phase. In the runtime phase, each sensor nodes takes 4 measurements of the signal emitted by the target, for a total of $m=144\ll n$ measurements. The considered model is the indoor model defined by the IEEE 802.15.4 standard, as reported in \cite[Equation 11]{fen09}. A measurement noise corresponding to an SNR of  $25$dB is added.
The online tracking is performed with O-IST, O-DR, and O-DISTA, which are run for 50 ms at each iteration, this time being sufficiently small such that, in the next measurement, the target is in the same cell or in adjacent cell. As to O-IST and O-DR, the data from sensors are processed in a centralized fusion center, while for O-DISTA the processing is performed in-network, with local communications based on the grid topology depicted in the third graph of Figure \ref{fig:s5}; specifically, each node can communicate with nodes at a maximum distance of 4.5 m.  We highlight that this topology is not $d$-regular, which goes beyond Assumption \ref{ass:dreg}.
\begin{figure*}
\centering
\includegraphics[width=0.49\columnwidth]{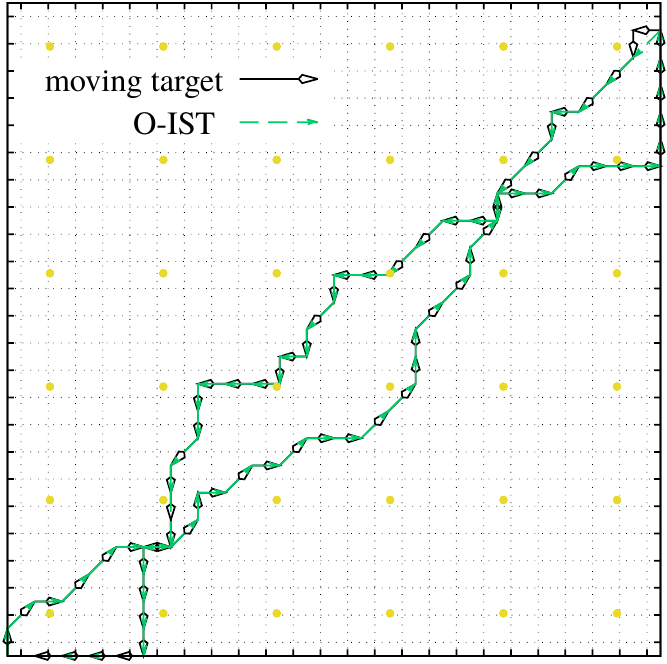}
\includegraphics[width=0.49\columnwidth]{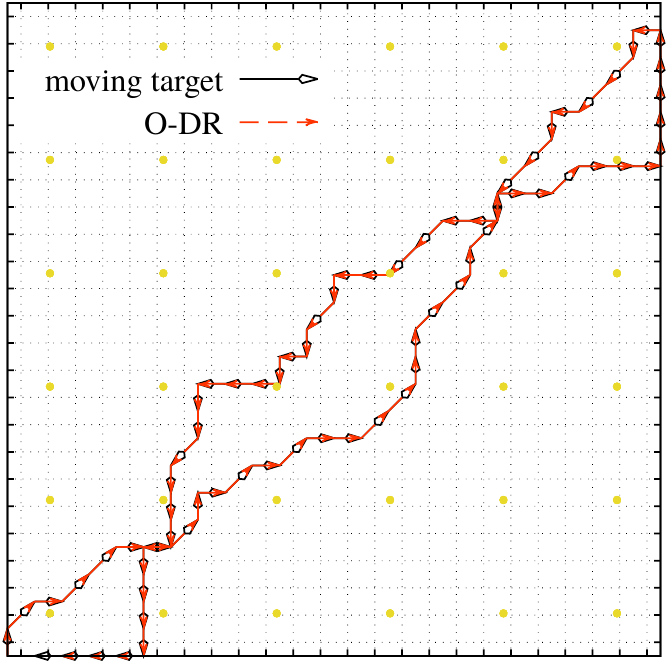}$~$
\includegraphics[width=0.49\columnwidth]{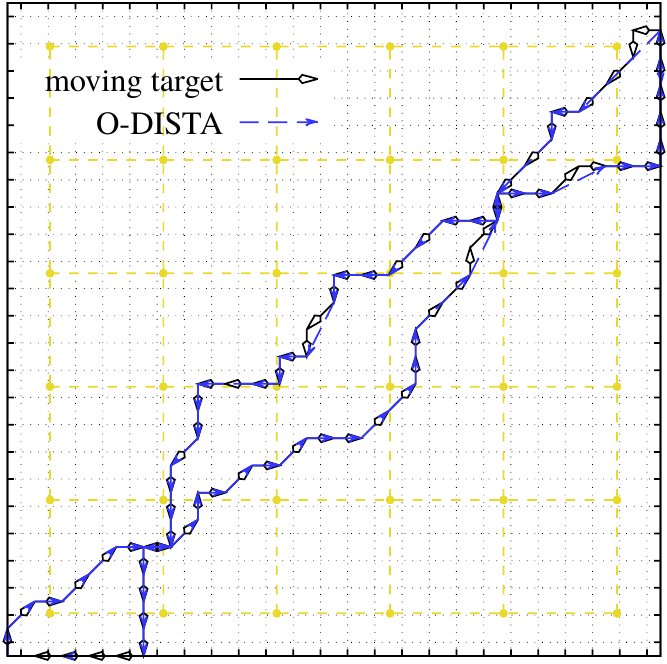}
\caption{Indoor tracking of a moving target from RSS compressed measurements, in a $25\times 25~\text{m}^2$ area. The path of the target and the corresponding online estimations are depicted. The yellow points denote the sensor nodes.}
	\label{fig:s5}
\end{figure*}
\begin{figure}
\centering
\includegraphics[width=0.49\columnwidth]{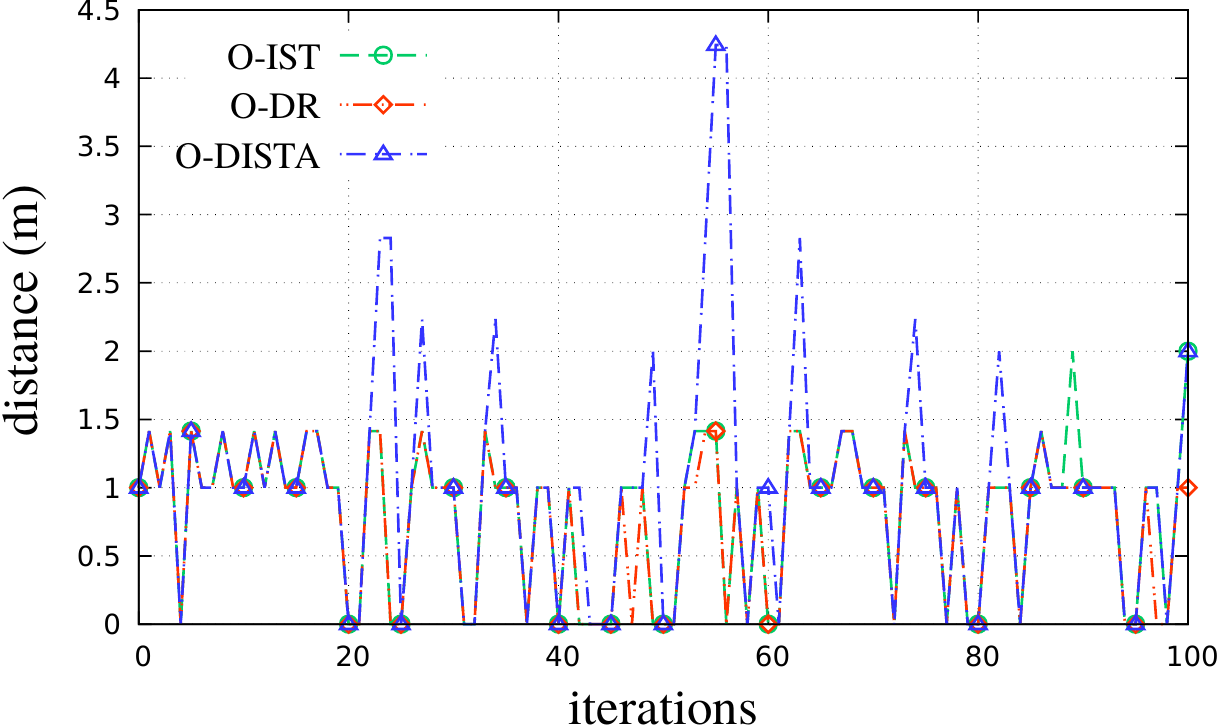}
\includegraphics[width=0.49\columnwidth]{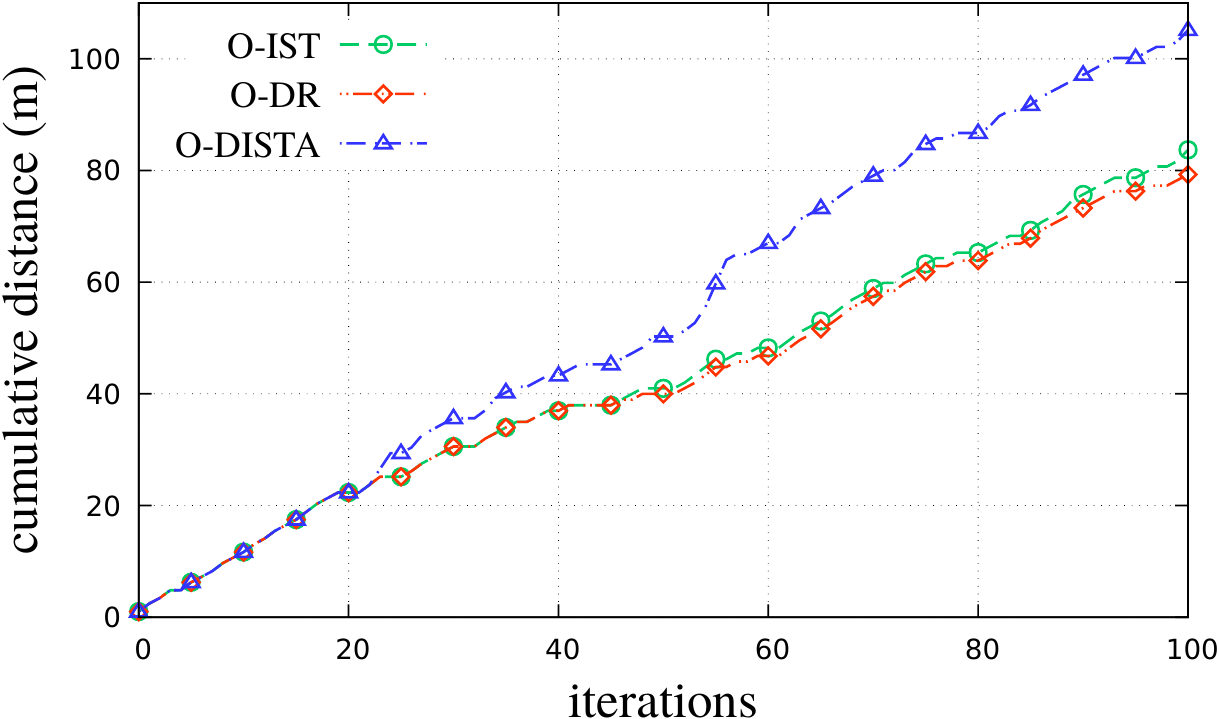}
\caption{Indoor tracking of a moving target from RSS compressed measurements: instantaneous distance and cumulative distance.}
	\label{fig:s6}
\end{figure}
In Figure \ref{fig:s5}, we show the path of the target in the $25\times 25~\text{m}^2$ area, and the corresponding tracking. We see that the three algorithms are substantially able to track the path. O-DR is the most precise and responsive, while O-DISTA is a bit less accurate. The performance is better visible in Figure \ref{fig:s6}, where the distance $\|x_t-\xtrue_t\|_2$ and the cumulative distance $\sum_{i=1}^t\|x_i-\xtrue_i\|_2$ are shown. A distance of $1$ m or of $\sqrt{2}$ m is natural when the target is moving, since the delay for processing in envisaged: such distances mean that the target has moved in an adjacent (horizontal/vertical or diagonal) cell. The distance may be null when the target stops in a cell. On the basis of this observation, O-DR is optimal, while O-DIST and O-DISTA present few more delays and missed corners.

Finally, we remark that, in this experiment the dictionary $A$ is constant, then in Assumption \ref{ass:relaxed_boundedness}, $Q_t$ is constant, which means that $\xtrue_t$ is not required to be bounded. In practice, this means that the moving target $\xtrue_t$ is not required to move within a fixed area to match the theoretical features of O-DISTA. On the other hand, if the area is not priorly fixed, a system of moving sensors should be provided.
\section{Conclusions}
In this work, we develop and analyze novel centralized and distributed strategies for sparse time-varying optimization. Specifically, we consider quadratic, strongly convex, optimization problems with $\ell_1$ regularization. In this setting, we provide a rigorous analysis in terms of dynamic regret. Furthermore,   numerical experiments on compressed system identification and indoor moving target tracking are presented. Future work will be devoted to extend the analysis to larger classed of problems. 
\bibliographystyle{plain}
\bibliography{refs}

\begin{thebibliography}{10}

\bibitem{agu17}
R.~P. {Aguilera}, G.~{Urrutia}, R.~A. {Delgado}, D.~{Dolz}, and J.~C.
  {Ag\"{u}ero}.
\newblock Quadratic model predictive control including input cardinality
  constraints.
\newblock {\em IEEE Trans. Autom. Control}, 62(6):3068--3075, 2017.

\bibitem{akb19}
M.~{Akbari}, B.~{Gharesifard}, and T.~{Linder}.
\newblock Individual regret bounds for the distributed online alternating
  direction method of multipliers.
\newblock {\em IEEE Trans. Autom. Control}, 64(4):1746--1752, 2019.

\bibitem{bal15}
A.~Balavoine, C.~J. Rozell, and J.~Romberg.
\newblock Discrete and continuous-time soft-thresholding for dynamic signal
  recovery.
\newblock {\em IEEE Trans. Signal Process.}, 63(12):3165--3176, 2015.

\bibitem{bec09}
A.~Beck and M.~Teboulle.
\newblock A fast iterative shrinkage-thresholding algorithm for linear inverse
  problems.
\newblock {\em SIAM J. Imaging Sci.}, 2(1):183--202, 2009.

\bibitem{bed18}
A.~S. Bedi, P.~Sarma, and K.~Rajawat.
\newblock Tracking moving agents via inexact online gradient descent algorithm.
\newblock {\em IEEE J. Sel. Topics Signal Process.}, 12(1):202--217, 2018.

\bibitem{ber18}
A.~Bernstein and E.~Dall'Anese.
\newblock Asynchronous and distributed tracking of time-varying fixed points.
\newblock In {\em Proc. IEEE Conf. Decis. Control (CDC)}, pages 791--798, 2018.

\bibitem{boy10}
S.~Boyd, N.~Parikh, E.~Chu, B.~Peleato, and J.~Eckstein.
\newblock Distributed optimization and statistical learning via the alternating
  direction method of multipliers.
\newblock {\em Found. Trends Mach. Learn.}, 3(1):1 -- 122, 2010.

\bibitem{bra16}
J.~M. Bravo, A.~Suarez, M.~Vasallo, and T.~Alamo.
\newblock Slide window bounded-error time-varying systems identification.
\newblock {\em IEEE Trans. Autom. Control}, 61(8):2282--2287, 2016.

\bibitem{cao19}
X.~{Cao} and K.~J.~R. {Liu}.
\newblock Dynamic sharing through the {ADMM}.
\newblock {\em IEEE Trans. Autom. Control}, page Early Access, 2019.

\bibitem{cha16}
A.~S. Charles, A.~Balavoine, and C.~J. Rozell.
\newblock Dynamic filtering of time-varying sparse signals via $\ell _1$
  minimization.
\newblock {\em IEEE Trans. Signal Process.}, 64(21):5644--5656, 2016.

\bibitem{che17}
T.~Chen, Q.~Ling, and G.~B. Giannakis.
\newblock An online convex optimization approach to proactive network resource
  allocation.
\newblock {\em IEEE Trans. Signal Process.}, 65(24):6350--6364, 2017.

\bibitem{dall19}
E.~Dall'Anese, A.~Simonetto, S.~Becker, and L.~Madden.
\newblock Optimization and learning with information streams: Time-varying
  algorithms and applications.
\newblock {\em arXiv 1910.08123}, 2019.

\bibitem{dar15}
D.~{Dardari}, P.~{Closas}, and P.~M. {Djuri\'{c}}.
\newblock Indoor tracking: Theory, methods, and technologies.
\newblock {\em IEEE Trans. Vehic. Tech.}, 64(4):1263--1278, 2015.

\bibitem{del16}
R.~R. {De Lucia}, S.~M. {Fosson}, and E.~{Magli}.
\newblock Low-power distributed sparse recovery testbed on wireless sensor
  networks.
\newblock In {\em IEEE Intern. Workshop Multim. Signal Process. (MMSP)}, 2016.

\bibitem{dix19}
R.~{Dixit}, A.~S. {Bedi}, R.~{Tripathi}, and K.~{Rajawat}.
\newblock Online learning with inexact proximal online gradient descent
  algorithms.
\newblock {\em IEEE Trans. Signal Process.}, 67(5):1338--1352, 2019.

\bibitem{don06}
D.~L. Donoho.
\newblock Compressed sensing.
\newblock {\em IEEE Trans. Inf. Theory}, 52(4):1289--1306, 2006.

\bibitem{duc10}
J.~C. Duchi, S.~Shalev-Shwartz, Y.~Singer, and A.~Tewari.
\newblock Composite objective mirror descent.
\newblock In {\em Conf. Learn. Theory (COLT)}, pages 14--26, 2010.

\bibitem{fen12}
C.~{Feng}, W.~S.~A. {Au}, S.~{Valaee}, and Z.~{Tan}.
\newblock Received-signal-strength-based indoor positioning using compressive
  sensing.
\newblock {\em IEEE Trans. Mobile Comput.}, 11(12):1983--1993, 2012.

\bibitem{fen09}
C.~Feng, S.~Valaee, and Z.~Tan.
\newblock Multiple target localization using compressive sensing.
\newblock In {\em IEEE Global Telecommunications Conference (GLOBECOM)}, pages
  1--6, 2009.

\bibitem{fia18}
A.~Fiandrotti, S.~M. Fosson, C.~Ravazzi, and E.~Magli.
\newblock {GPU}-accelerated algorithms for compressed signals recovery with
  application to astronomical imagery deblurring.
\newblock {\em Int. J. Remote Sens.}, 39(7):2043--2065, 2018.

\bibitem{for10}
M.~Fornasier.
\newblock Numerical methods for sparse recovery.
\newblock In M.~Fornasier, editor, {\em Theoretical Foundations and Numerical
  Methods for Sparse Recovery}, pages 93--200. Radon Series Comp. Appl. Math.,
  de Gruyter, 2010.

\bibitem{fox18cdc}
S.~M. Fosson.
\newblock Online optimization in dynamic environments: a regret analysis for
  sparse problems.
\newblock In {\em Proc. IEEE Conf. Decis. Control (CDC)}, 2018.

\bibitem{fou13}
Simon Foucart and Holger Rauhut.
\newblock {\em A Mathematical Introduction to Compressive Sensing}.
\newblock Springer, New York, 2013.

\bibitem{gal16}
Marco Gallieri.
\newblock {\em Lasso-MPC - Predictive Control with $\ell_1$-Regularised Least
  Squares}.
\newblock Springer Publishing Company, Incorporated, 1st edition, 2016.

\bibitem{gis17}
P.~Giselsson and S.~Boyd.
\newblock Linear convergence and metric selection for {D}ouglas-{R}achford
  splitting and {ADMM}.
\newblock {\em IEEE Trans. Autom. Control}, (62), 2017.

\bibitem{hal15}
E.~C. Hall and R.~M. Willett.
\newblock Online convex optimization in dynamic environments.
\newblock {\em IEEE J. Sel. Topics Signal Process.}, 9(4):647--662, 2015.

\bibitem{haz07}
Elad Hazan, Amit Agarwal, and Satyen Kale.
\newblock Logarithmic regret algorithms for online convex optimization.
\newblock {\em Machine Learning}, 69(2-3):169--192, 2007.

\bibitem{hos16}
S.~Hosseini, A.~Chapman, and M.~Mesbahi.
\newblock Online distributed convex optimization on dynamic networks.
\newblock {\em IEEE Trans. Autom. Control}, 61(11):3545--3550, 2016.

\bibitem{jer14}
J.~L. Jerez, P.~J. Goulart, S.~Richter, G.~A. Constantinides, E.~C. Kerrigan,
  and M.~Morari.
\newblock Embedded online optimization for model predictive control at
  megahertz rates.
\newblock {\em IEEE Trans. Autom. Control}, 59(12):3238--3251, 2014.

\bibitem{kha17}
A.~{Khalajmehrabadi}, N.~{Gatsis}, D.~J. {Pack}, and D.~{Akopian}.
\newblock A joint indoor {WLAN} localization and outlier detection scheme using
  {LASSO} and elastic-net optimization techniques.
\newblock {\em IEEE Trans. Mobile Comput.}, 16(8):2079--2092, 2017.

\bibitem{lee00}
Jay~H. Lee, Kwang~S. Lee, and Won~C. Kim.
\newblock Model-based iterative learning control with a quadratic criterion for
  time-varying linear systems.
\newblock {\em Automatica}, 36(5):641 -- 657, 2000.

\bibitem{lee18}
S.~Lee, A.~Nedi\'{c}, and M.~Raginsky.
\newblock Coordinate dual averaging for decentralized online optimization with
  nonseparable global objectives.
\newblock {\em IEEE Trans. Control Netw. Syst.}, 5(1):34--44, 2018.

\bibitem{li11}
Y.~Li, H.~l.~Wei, and S.~A. Billings.
\newblock Identification of time-varying systems using multi-wavelet basis
  functions.
\newblock {\em IEEE Trans. Control Syst. Technol.}, 19(3):656--663, 2011.

\bibitem{linrib14}
Q.~Ling and A.~Ribeiro.
\newblock Decentralized dynamic optimization through the alternating direction
  method of multipliers.
\newblock {\em IEEE Trans. Signal Process.}, 62(5):1185--1197, 2014.

\bibitem{lio79}
P.~L. Lions and B.~Mercier.
\newblock Splitting algorithms for the sum of two nonlinear operators.
\newblock {\em SIAM J. Numer. Anal.}, (6):964--979, 1979.

\bibitem{lu19}
K.~{Lu}, G.~{Jing}, and L.~{Wang}.
\newblock Online distributed optimization with strongly pseudoconvex-sum cost
  functions.
\newblock {\em IEEE Trans Autom. Control}, pages 1--1, 2019.

\bibitem{mar18}
M.~{Maros} and J.~{Jald\'{e}n}.
\newblock {ADMM} for distributed dynamic beamforming.
\newblock {\em IEEE Trans. Signal Inf. Process. Netw.}, 4(2):220--235, 2018.

\bibitem{mar19}
Marie Maros and Joakim Jald\'{e}n.
\newblock On decentralized tracking with {ADMM} for problems with time-varying
  curvature.
\newblock {\em arXiv 1903.06492}, 2019.

\bibitem{mata15}
J.~Matamoros, S.~M. Fosson, E.~Magli, and C.~Ant\'{o}n-Haro.
\newblock Distributed {ADMM} for in-network reconstruction of sparse signals
  with innovations.
\newblock {\em IEEE Trans. Signal Inf. Process. Netw.}, 1(4):225--234, 2015.

\bibitem{mok16}
Aryan Mokhtari, Shahin Shahrampour, Ali Jadbabaie, and Alejandro Ribeiro.
\newblock Online optimization in dynamic environments: Improved regret rates
  for strongly convex problems.
\newblock In {\em Proc. IEEE Conf. Decis. Control (CDC)}, 2016.

\bibitem{mot15}
J.~Mota, N.~Deligiannis, A.~C. Sankaranarayanan, V.~Cevher, and M.~Rodrigues.
\newblock Dynamic sparse state estimation using $\ell_1-\ell_1$ minimization:
  Adaptive-rate measurement bounds, algorithms and applications.
\newblock In {\em IEEE Int. Conf. Acoust. Speech Signal Process. (ICASSP)},
  pages 3332--3336, 2015.

\bibitem{mou19}
Walaa~M. Moursi and Yuriy Zinchenko.
\newblock A note on the equivalence of operator splitting methods.
\newblock In Heinz~H. Bauschke, Regina~S. Burachik, and D.~Russell Luke,
  editors, {\em Splitting Algorithms, Modern Operator Theory, and
  Applications}, pages 331--349. Springer International Publishing, 2019.

\bibitem{oom17}
Tom Oomen and Cristian~R. Rojas.
\newblock Sparse iterative learning control with application to a wafer stage:
  Achieving performance, resource efficiency, and task flexibility.
\newblock {\em Mechatronics}, 47:134--147, 2017.

\bibitem{pra16}
G.~Prando, D.~Romeres, and A.~Chiuso.
\newblock Online identification of time-varying systems: A {B}ayesian approach.
\newblock In {\em Proc. IEEE Conf. Decis. Control (CDC)}, pages 3775--3780,
  2016.

\bibitem{rah17}
S.~Rahili and W.~Ren.
\newblock Distributed continuous-time convex optimization with time-varying
  cost functions.
\newblock {\em IEEE Trans. Autom. Control}, 62(4):1590--1605, 2017.

\bibitem{rfm15}
C.~Ravazzi, S.~M. Fosson, and E.~Magli.
\newblock Distributed iterative thresholding for $\ell_0$/$\ell_1$-regularized
  linear inverse problems.
\newblock {\em IEEE Trans. Inf. Theory}, 61(4):2081--2100, 2015.

\bibitem{rav15}
C.~Ravazzi, S.~M. Fosson, and E.~Magli.
\newblock Randomized algorithms for distributed nonlinear optimization under
  sparsity constraints.
\newblock {\em IEEE Trans. Signal Process.}, 64(6):1420--1434, 2015.

\bibitem{san11}
B.~M. Sanandaji, T.~L. Vincent, M.~B. Wakin, and R.~T\'{o}th.
\newblock Compressive system identification of {LTI} and {LTV ARX} models.
\newblock In {\em Proc. IEEE Conf. Decis. Control (CDC)}, pages 791--798, 2011.

\bibitem{sha18}
S.~Shahrampour and A.~Jadbabaie.
\newblock Distributed online optimization in dynamic environments using mirror
  descent.
\newblock {\em IEEE Trans. Autom. Control}, 63(3):714--725, 2018.

\bibitem{sha12book}
Shai Shalev-Shwartz et~al.
\newblock Online learning and online convex optimization.
\newblock {\em Found. Trends Mach. Learn.}, 4(2):107--194, 2012.

\bibitem{sim17pre}
A.~Simonetto.
\newblock Time-varying convex optimization via time-varying averaged operators.
\newblock {\em arXiv 1704.0733}, 2017.

\bibitem{sim18acc}
A.~Simonetto.
\newblock Prediction-correction dual ascent for time-varying convex programs.
\newblock In {\em Proc. Amer. Control Conf. (ACC)}, pages 4508--4513, 2018.

\bibitem{sim19}
A.~Simonetto.
\newblock Dual prediction-correction methods for linearly constrained
  time-varying convex programs.
\newblock {\em IEEE Trans. Autom. Control}, 64(8):3355--3361, 2019.

\bibitem{sim17tsp}
A.~Simonetto and E.~Dall'Anese.
\newblock Prediction-correction algorithms for time-varying constrained
  optimization.
\newblock {\em IEEE Trans. Signal Process.}, 65(20):5481--5494, 2017.

\bibitem{sim17}
A.~Simonetto, A.~Koppel, A.~Mokhtari, G.~Leus, and A.~Ribeiro.
\newblock Decentralized prediction-correction methods for networked
  time-varying convex optimization.
\newblock {\em IEEE Trans. Autom. Control}, 62(11):5724--5738, 2017.

\bibitem{sim16TSP}
A.~Simonetto, A.~Mokhtari, A.~Koppel, G.~Leus, and A.~Ribeiro.
\newblock A class of prediction-correction methods for time-varying convex
  optimization.
\newblock {\em IEEE Trans. Signal Process.}, 64(17):4576--4591, 2016.

\bibitem{sop16}
P.~{Sopasakis}, N.~{Freris}, and P.~{Patrinos}.
\newblock Accelerated reconstruction of a compressively sampled data stream.
\newblock In {\em Proc. Eur. Signal Process. Conf. (EUSIPCO)}, pages
  1078--1082, 2016.

\bibitem{sun17}
C.~Sun, M.~Ye, and G.~Hu.
\newblock Distributed time-varying quadratic optimization for multiple agents
  under undirected graphs.
\newblock {\em IEEE Trans. Autom. Control}, 62(7):3687--3694, 2017.

\bibitem{suz13}
Taiji Suzuki.
\newblock Dual averaging and proximal gradient descent for online alternating
  direction multiplier method.
\newblock In {\em Proc. Int. Conf. Mach. Learn. (ICML)}, volume~28, pages
  392--400, 2013.

\bibitem{tib96}
R.~Tibshirani.
\newblock Regression shrinkage and selection via the lasso.
\newblock {\em Journal of the Royal Statistical Society, Series B},
  58:267--288, 1996.

\bibitem{tib13}
Ryan~J. Tibshirani.
\newblock {The {L}asso problem and uniqueness}.
\newblock {\em Elec. J. Stat.}, 7:1456--1490, 2013.

\bibitem{tot11}
R.~T\'{o}th, B.~M. Sanandaji, K.~Poolla, and T.~L. Vincent.
\newblock Compressive system identification in the linear time-invariant
  framework.
\newblock In {\em Proc. IEEE Conf. Decis. Control (CDC)}, pages 783--790, 2011.

\bibitem{vas16}
N.~Vaswani and J.~Zhan.
\newblock Recursive recovery of sparse signal sequences from compressive
  measurements: A review.
\newblock {\em IEEE Trans. Signal Process.}, 64(13):3523--3549, 2016.

\bibitem{wan12}
Huahua Wang and Arindam Banerjee.
\newblock Online alternating direction method.
\newblock In {\em Proc. Int. Conf. Mach. Learn. (ICML)}, pages 1119--1126,
  2012.

\bibitem{clairvoyant}
T.~Yang, L.~Zhang, R.~Jin, and J.~Yi.
\newblock Tracking slowly moving clairvoyant: Optimal dynamic regret of online
  learning with true and noisy gradient.
\newblock In {\em Proc. Int. Conf. Mach. Learn. (ICML)}, pages 449--457, 2016.

\bibitem{zac12}
D.~Zachariah, S.~Chatterjee, and M.~Jansson.
\newblock Dynamic iterative pursuit.
\newblock {\em IEEE Trans. Signal Process.}, 60(9):4967--4972, 2012.

\bibitem{zha18}
L.~Zhang, T.~Yang, R.~Jin, and Z.-H. Zhou.
\newblock Dynamic regret of strongly adaptive methods.
\newblock In {\em Proc. Int. Conf. Mach. Learn. (ICML)}, pages 9372--9381,
  2018.

\bibitem{zha17}
L.~Zhang, T.~Yang, J.~Yi, R.~Jin, and Z.-H. Zhou.
\newblock Improved dynamic regret for non-degenerate functions.
\newblock In {\em Proc. Conf. Neural Inf. Process. Syst. (NIPS)}, pages
  733--742, 2017.

\bibitem{dall17}
Y.~{Zhang}, E.~{Dall'Anese}, and M.~{Hong}.
\newblock Dynamic {ADMM} for real-time optimal power flow.
\newblock In {\em Proc. IEEE Global Conf. Signal Inf. Process. (GlobalSIP)},
  pages 1085--1089, 2017.

\bibitem{zin13}
J.~Ziniel and P.~Schniter.
\newblock Dynamic compressive sensing of time-varying signals via approximate
  message passing.
\newblock {\em IEEE Trans. Signal Process.}, 61(21):5270--5284, 2013.

\bibitem{zin03}
M.~Zinkevich.
\newblock Online convex programming and generalized infinitesimal gradient
  ascent.
\newblock In {\em Proc. Int. Conf. Mach. Learn. (ICML)}, pages 928--936, 2003.

\bibitem{zou05}
H.~Zou and T.~Hastie.
\newblock Regularization and variable selection via the elastic net.
\newblock {\em J. Royal Stat. Soc. B}, 67(2):301--320, 2005.

\end{thebibliography}

\end{document}